\documentclass[12pt, reqno]{amsart}
\usepackage{amsmath, amsthm, amscd, amsfonts, amssymb}
\usepackage{bm}
\usepackage[all,cmtip]{xy}
\newtheorem{theorem}{Theorem}[section]

\newtheorem{corollary}[theorem]{Corollary}
\theoremstyle{definition}

\newtheorem{problem}[theorem]{Problem}

\numberwithin{equation}{section}
\newtheorem{Question}[theorem]{Question}

\newcommand{\overbar}[1]{\mkern 1.5mu\overline{\mkern-1.5mu#1\mkern-1.5mu}\mkern 1.5mu}
\newcommand\bs[1]{{\boldsymbol{#1}}}
\newcommand\mc[1]{{\mathcal{#1}}}
\newcommand\mb[1]{{\mathbb{#1}}}

\begin{document}

\setcounter{page}{1}

\title{A compendium of research in operator algebras and operator theory}

\author[Sarkar]{Jaydeb Sarkar}
\address{Indian Statistical Institute, Statistics and Mathematics Unit, 8th Mile, Mysore Road, Bangalore, 560059, India}\email{jaydeb@gmail.com, jay@isibang.ac.in}

\begin{abstract}
This chapter surveys the advances of the past decade arising from the contributions of Indian mathematicians in the broad areas of operator algebras and operator theory. It brings together the work of twenty mathematicians and their collaborators, each writing from the perspective of their respective research fields and beyond. Several problems highlighted here are expected to shape the future development of the subject at a global level.
\end{abstract}

\maketitle

\tableofcontents

\section*{Foreword}

Operator algebras and operator theory are subfields of analysis, more specifically of functional analysis. Broadly speaking, these areas concern linear operators, as well as the algebras, subalgebras, and vector spaces formed by such operators. The origins of the subject lie in the study of matrices and matrix spaces. It is also important to note that the nature of the underlying spaces, even at the level of finite-dimensional vector spaces, plays a significant role in the broader theory. This perspective is particularly apparent in the theory of Banach spaces.

The core of these areas touches upon nearly all major branches of mathematics, including harmonic analysis, partial and ordinary differential equations, complex analysis (both one and several variables), geometry, topology, and even number theory. The broader field also connects with control theory, branches of engineering, and recent developments in quantum computation and machine learning. The central aim of this area of linear analysis is to investigate the structure of linear operators as well as the spaces generated by such operators, much as one does for matrices. It has produced a number of results and connections of broad mathematical interest. Notable examples include the Atiyah–Singer index theorem; Brown–Douglas–Fillmore theory on the essential spectrum and index theory; the solution of the Bieberbach conjecture; the algebraic and topological approaches to the classification of $C^*$-algebras; the solution of the Kadison–Singer problem; Taylor’s Koszul complex approach to the spectrum of tuples of operators; the profound theory of the geometry of Banach spaces; wavelet theory; and the role of algebraic varieties and complex geometry in understanding the structure of operators in several variables, among other developments.

The theory of linear analysis, or more specifically matrix theory, was formally initiated through the work of Grassmann and Sylvester. However, it was Stefan Banach and John von Neumann who laid much of the foundation for the modern theory of linear operators on Hilbert and Banach spaces. Not long after this period, the culture of linear analysis in India also flourished and soon became a subject of great interest to many Indian mathematicians. The subject was then led by a number of mathematicians working in India including K.R. Parthasarathy, M.G. Nadkarni, Ashoke Kumar Roy, K.B. Sinha, Rajendra Bhatia, V.S. Sunder, Gadadhar Misra, Subhash J. Bhatt, B. V. Limaye, and Ameer Athavale, among others.

The baton from these senior mathematicians has been carried forward by a much larger and still-growing community working in operator algebras and operator theory in India today. In particular, the past decade has witnessed remarkable results from Indian mathematicians, including a surge of contributions from younger researchers.

It is both timely and important to report on the achievements of Indian mathematicians in these areas in the past decade and to highlight future directions and possibilities. This chapter consists of 19 contributions, along with a special extended section of open problems contributed by Gadadhar Misra. Together, these contributions cover a broad spectrum of mathematics within operator algebras, operator theory, and related topics. All the contributors are well-established researchers who have been active in the field for a considerable time. Collectively, their work showcases the present strength and the emerging future leadership of Indian mathematics on the global stage.

It is important to emphasize that our role in this chapter is primarily to bring together all the contributors on a common platform and to organize their work in a coherent manner. That said, the full scientific credit for the content of this chapter belongs to the authors of the 20 individual sections. The names, affiliations, and contact details of the contributors to each section are provided immediately after its title.

\section{Some problems in operator theory and related topics}

\noindent \textit{By Gadadhar Misra [gm@isibang.ac.in] from ISI Bangalore and IIT Gandhinagar.}

\bigskip

{\bf Introduction.} The development of multi-variate operator theory has  paralleled that of function theory in several complex variables except that the non-commuting co-ordinates pose an entirely new set of challenges. The introduction of methods of commutative algebra for studying problems in this area by viewing a pair $(\mathcal H, \boldsymbol T)$, where $\boldsymbol T:=(T_1, \ldots, T_n)$  is a commuting $n$-tuple of operators on the Hilbert space $\mathcal H$,
as a Hilbert module via the multiplication $(p,h) \mapsto p(\boldsymbol T) h$, $p\in \mathbb C[z_1, \ldots , z_n],\,h\in\mathcal H$,
was the beginning of a systematic new development. This perspective is powerful because it allows the application of tools from commutative algebra. On the other hand, these techniques don't apply directly because of the   continuity assumption of the module multiplication either in just the first variable, or  in both the variables. A choice is made depending on the problem at hand.

A second ingredient has been the action of a group locally compact second countable topological group $G$ on the module $(\mathcal H, \boldsymbol T)$.
Assume that $X$ is a nonempty subset $\mathbb{C}^d$ and that the action of $g\in G$ is holomorphic in some open neighbourhood of $\overline{X}$. If $\boldsymbol{T}$ is a $d$- tuple and $X$ is the Taylor joint spectrum of $\sigma(\boldsymbol{T})$, then  $g\cdot \boldsymbol T$ is defined by the usual holomorphic functional calculus. The classification of imprimitivity systems $\big (X, G, (\mathcal H,\boldsymbol T)\,\big ):=  \big \{\boldsymbol T: U_g^* p(\boldsymbol T) U_g = (g\cdot p) (\boldsymbol T)\},$ defined as the set of $d$- tuples $\boldsymbol{T}$ satisfying $U^{*}{g}p(\boldsymbol{T})U{g}=(g\cdot p)(\boldsymbol{T})$ for all $g\in G, p\in \mathbb{C}[\boldsymbol{z}]$; where
$U_g$ is a unitary representation of the group $G$ and $g\cdot p = p\circ g^{-1}$, is an important problem. One may study such imprimitivities by choosing a group $G$ acting on $X$. This choice might vary from a transitive action to an action, where the number of orbits of $G$ in $X$ are quite large. For a typical example, consider $X=\overbar{\mathbb D}^n$, the $n$-fold product of the unit closed  disc, and the group $G$ to be either the bi-holomorphic automorphism group of $\overbar{\mathbb D}^n$ or simply the permutation group $\mathfrak S_n$ acting on $\overbar{\mathbb D}^n$.

The study of these problems are greatly facilitated by the analysis of a class holomorphic curves in the Grassmannian $\mbox{\rm Gr}(\mathcal H,n)$ of rank $n$ in some separable complex Hilbert space $\mathcal H$. These holomorphic curves arise from a  class of operators acting on some Hilbert space $\mathcal H$  introduced in the paper \cite{CD}.  The operators in  this  class  possess  an open set $\Omega\subseteq \mathbb C$ of eigenvalues of (constant) multiplicity $k$ and characterized by the existence of a holomorphic map $\gamma:\Omega \to \mathcal H$ such that $\boldsymbol \gamma(w):=(\gamma_1(w), \ldots, \gamma_k(w))$,
$T\gamma_i(w) = w\gamma_i(w),\: 1\leq i \leq k,\: w\in \Omega.$
For $k=1$, one of the main features of the operator $T$ in this class is that the curvature
$
\mathcal K_T(w):= -\partial \bar{\partial}\log \|\gamma_T(w)\|^2
$
of the holomorphic Hermitian line bundle $E_T$ determined by the  holomorphic map $\gamma_T$ equipped with the Hermitian structure $\|\gamma_T(w)\|^2$ is a complete unitary invariant for the operator $T$.

A complete set of invariants, originally  obtained by Cowen and Douglas for any $k \geq  1$, have been refined to provide a tractable set of invariants for large class of Cowen-Douglas  operators in $B_k(\Omega)$, $\Omega \subseteq \mb{C}$, \cite{JJKM, JJM}.
In these papers, after imposing a mild condition on the Cowen-Douglas bundles, it is shown that the curvature together with the second fundamental form serves as a complete set of invariants.


{\bf Homogeneous operators.} It is easy to see that if $T$ is a contraction in the Cowen-Douglas class of the unit disc $\mathbb D$, then $\mathcal K_T(w) \leq \mathcal K_{S^*}(w),$ where $S^*$ is the backward unilateral shift acting on $\ell^2$. Choosing a holomorphic frame $\gamma_{S^*},$ say $\gamma_{S^*}(w)=(1,w,w^2, \ldots )$, it follows that
$$\|\gamma_{S^*}(w)\|^2 = (1-|w|^2)^{-1},\,\, \mathcal K_{S^*}(w) =-(1-|w|^2)^{-2},\,\, w\in \mathbb D,$$ with respect to the frame $\gamma$. Thus the operator $S^*$ is an extremal operator in the class of all contractive Cowen-Douglas operator.  R. G. Douglas asked if the curvature $\mathcal K_T$ of a contraction $T$ achieves equality in this inequality even at just one point, then does it follow that $T$ must be unitarily equivalent to $S^*$?  It is easy to see that the answer is ``no'', in general. However, if $T$ is homogeneous, namely, $U_\varphi^* T U_\varphi = \varphi(T)$ for each bi-holomorphic automorphism $\varphi$ of the unit disc and some unitary $U_\varphi$, then the answer is ``yes''. Of course, it is then natural to ask what are all the homogeneous operators. 

Assume there exists a projective unitary representation $\sigma$ of M\"{o}b
such that
$\varphi(T)= \sigma(\varphi)^\star T \sigma(\varphi)$ for all  for all M\"{o}bius transformations $\varphi$.  If this is the case, the operator $T$ is homogeneous and we say that $\sigma(\varphi)$  is associated with the operator $T$.
A M\"{o}bius equivariant version of the  Sz.-Nagy--Foias model theory for completely non-unitary (cnu) contractions is developed in [67].  As an application, we prove that if T is a cnu contraction with associated (projective
unitary) representation $\sigma$, then there is a unique projective unitary
representation $\hat{\sigma}$, extending $\sigma$, associated with the minimal unitary dilation  of $T$. The representation $\hat{\sigma}$, which extends $\sigma$, is given in terms of $\sigma$ by the formula
$$
\hat{\sigma} = (\pi \otimes D_1^+) \oplus \sigma \oplus (\pi_\star \otimes D_1^-),
$$
where $D_1^\pm$  are the two Discrete series representations (one holomorphic and the other anti-holomorphic) living on the Hardy space $H^2(\mathbb D)$, and $\pi, \pi_\star$  are representations of M\"{o}b living on the two defect spaces of $T$  and defined explicitly in terms of $\sigma$ and $T$.
Moreover, a cnu contraction $T$ has an associated representation
if and only if its Sz.-Nagy--Foias characteristic function $\theta_T$ has the product form
$\theta_T(z) = \pi_\star(\varphi_z)^* \theta_T(0) \pi(\varphi_z),$ $z\in \mathbb D$, where $\varphi_z$ is the involution in M\"{o}b mapping $z$ to $0.$
We obtain a concrete realization of this product formula
for a large subclass of homogeneous cnu contractions from the Cowen-Douglas class; these are the holomorphic imprimitivity systems among the homogeneous contractions.

\begin{problem}
Find all the homogeneous contractions using the product formula for the Sz.-Nagy -- Foias characteristic function.
\end{problem}

Given a pair of positive real numbers $\alpha, \beta$ and a sesqui-analytic function $K$ on a bounded domain $\Omega \subseteq \mathbb C^m$, we investigate the properties of the real-analytic function
$$K^{(\alpha,\beta)}(\boldsymbol z, \boldsymbol z):=\big (\!\big ( \tfrac{\partial^2}{\partial \overbar{z}_j \partial z_i } K^{(\alpha+\beta )} \log K (\boldsymbol z,\boldsymbol z) \big )\!\big )_{1\leq i,j \leq m},\,\, \boldsymbol z\in \Omega,$$
taking values in $m \times m$ matrices.  The kernel $K^{(\alpha,\beta)}$ is non-negative definite whenever $K^\alpha$ and $K^\beta$ are non-negative definite.  In this case, a realization of the Hilbert module determined by the kernel $K^{(\alpha,\beta)}$ is obtained. Let $\mathcal M_i$, $i=1,2$, be two Hilbert modules over the polynomial ring $\mathbb C[z_1, \ldots , z_m]$.  The tensor product $\mathcal M_1\otimes \mathcal M_2$ is clearly a module over the ring $\mathbb C[z_1, \ldots , z_{2m}]$.  This module multiplication  restricts to $\mathbb C[z_1, \ldots , z_m]$ via the diagonal map $\boldsymbol z \mapsto (\boldsymbol z, \boldsymbol z)$, $\boldsymbol z \in \Omega$. Now, a natural decomposition of the tensor product $\mathcal M_1\otimes \mathcal M_2$ very similar to the Clebsch-Gordon decomposition of the tensor product of two irreducible unitary representations occurs relative to the multiplication restricted to
$\mathbb C[z_1, \ldots , z_m]$. Two of the initial pieces in this decomposition have been identified in \cite{GM}.
The first of these is simply the restriction of $\mathcal M_1\otimes \mathcal M_2$ to the diagonal set $\bigtriangleup:=\{ (\boldsymbol z, \boldsymbol z): \boldsymbol z\in \Omega\}$, while the second piece in the decomposition is the module determined by the kernel $K^{(\alpha,\beta)}$. Moreover, if $\Omega$ is a bounded symmetric domain, then $K^{(\alpha,\beta)}$ is covariant whenever $K$ is covariant under the action of the bi-holomorphic automorphism group of $\Omega$.
\begin{problem}
Find all the components in the decomposition of the tensor product of two Hilbert modules. This would be analogous to the Clebsch - Gorden decomposition for tensor product of two unitary representations.
\end{problem}


{\bf Holomorphic imprimitivity.} Let $(\mathcal H, \boldsymbol T)$ be a Hilbert module over a ring $\mc{R}$. Suppose that a group $G$ acts on $\mc{H}$, that is, there exists a map $U: G \times \mc{H} \to \mc{H}$ defined by $U(g,f) = U_g f$ for some unitary homomorphism $U_g$, $g\in G$, $f\in \mc{H}$.  Assume also that the group $G$ acts on the ring $\mc{R}$ and let us denote this action by $(g, \phi) \to g\cdot \phi$, $g\in G$ and $\phi \in \mc{R}$.  The imprimitivity condition is the requirement
\[ U_g^* m_\phi(f)  U_g = m_{g\cdot \phi}(f),\,\, g\in G, h \in \mc{H}, \phi \in \mc{R}. \]
Specializing to the case where the ring $\mc{R}$ is a commutative $C^*$ algebra, that is, $\mc{R)} = C(X)$ for some locally compact Hausdroff space $X$, where $X$ is a $G$ - space gives the Mackey imprimitivity.  In this case, the module $(\mathcal H, \boldsymbol T)$ can be identified with $L^2(X, d\mu)$ and $\bs{T}$ is the commuting tuple of multiplication operators on $L^2(X, d\mu)$, where $d\mu$ is a measure on $X$ quasi-invariant under the action of the group $G$. Thus the commuting tuple $\bs{T}$ is a commuting tuple of normal operators. The condition of imprimitivity is the requirement that $\bs{T}$ is homogeneous, that is,
$$U_g^*  \boldsymbol T U_g = g\cdot \boldsymbol T,\,\, g\in G,$$
where the continuous action of $G$ on $X$ is nothing but a continuous function on $X$ and the $g\cdot \boldsymbol T$ is defined via the usual continuous functional calculus for a commuting tuple of normal operators. The imprimitivity has an apparent stronger form. The action of $G$ on $X$ induces a $*$ - homomorphism of $C(X)$ given by the formula $\phi_g: \phi \to \phi\circ g,\,\phi\in C(X)$. Now, the imprivity relation can be written in the form
$$U_g^*  \phi(\boldsymbol T) U_g = \phi_g(\boldsymbol T),\,\, g\in G,\, \phi \in C(X).$$
Mackey's imprimitivity theorem says that these imprimitivities are in one to one correspondence with induced representations of the subgroup $H$, where $H$ is determined by realizing $X$ in the form $X:= G/H$. An imprimitivity system is called irreducible if there is no closed subspace $\mc{M} \subseteq \mc{H}$ such that $\mc{M}$ is invariant under the representation $U_g$ of the group $G$ and the $*$ - representation $\phi \to \phi(\bs{T})$.  One of the early results of Mackey is that irreducible represenations of the subgroup $H$ corresponds to irreducible multiplier representations of the group $G$. There is a vast literature on this subject. However, there has been no attempt to understand the subspaces $\mc{M}$ invariant under both $U_g$ and $M_\phi$, $g\in G$ and $\phi \in C(X)$.  The restriction of an imprimitivity to such an invariant subspace corresponds to homogeneous subnormal operators. In joint work with A. Kor\'{a}nyi, we are in the process of describing these invaraint subpaces in one of the simplest cases, namely, where $X = \mathbb D$ is the open unit disc and $G$ is the M\"{o}bius group of biholomorphic automorphisms of the disc $\mathbb D$. These invariant subspaces, typically, consist of holomorphic functions defined on $X$, hence we choose to call the restriction of an imprimitivity (in the sense of Mackey) to an invariant subspace ``holomorphic imprimitivity".

Consider for example, the Hilbert space  $L^{(\lambda)}$, $\lambda \in \mathbb R$, of measurable  (complex-valued) functions on the open unit disc $\mathbb D$,  such that
$$
\|f\|^2_\lambda = \int_{\mathbb D} (1-|z|^2)^{\lambda-2} |f(z)|^2 dA_z < \infty.
$$
($dA_z$ stands for $dxdy$, $z=x+iy$.)  Let $G$ be the simply connected covering group of $SU(1,1)$; $g(z)$ and $g^\prime(z)=\tfrac{\partial(g(z))}{\partial z}$ still make sense for $g\in G$.
The group $G$ acts on $L^{(\lambda)}$ via a unitary representation ${U_{g}}$ defined by the multiplier $g^{\prime}(z)^{-\lambda/2}$. Explicitly (it is simpler to write down $U_{g^{-1}}$ than $U_{g}$),
$$ U_{g^{-1}} f(z) = g^\prime(z)^{\lambda/2} f (g(z)).$$
Unitarity follows from the identity $(1-|z|^2) |g^\prime(z)| = 1 - |g(z)|^2$.

The representation induced in the sense of Mackey by the character $-\tfrac{\lambda}{2}$ of $K$, the stabiilizer of $0$ in $G$, is $U_g$. This follows by observing that any $g\in K$ acts by a rotation $g(z) = e^{i\theta} z$, and then ${g^\prime(z)}^{-\lambda/2} = e^{-i\tfrac{\lambda}{2} \theta}$. The subspace $H^{(\lambda)}$ of holomorphic functions in $L^{(\lambda)}$ is invariant under $\{U_g: g\in G\}$.
This subspace $H^{(\lambda)}$ is not zero if and only if $\lambda > 1$; in this case, the restriction of $\{U_g:g\in G\}$ to $H^{(\lambda)}$ is the holomorphic discrete series. The subspace $H^{(\lambda)}$ is also an invariant subspace for the operator $M$ defined by $Mf(z) = z f(z)$.

\begin{problem} Is $H^{(\lambda)}$ the only subspace invariant under both $\{U_g: g\in G\}$ and  $M$? \end{problem}

The question of classifying the holomorphic imprimitivities on a  bounded symmetric domain $\Omega$ amounts to classification of commuting  tuples  of homogeneous of operators in the Cowen-Douglas class of $\Omega$. There is a one to one correspondence  between these and the holomomorphic homogeneous vector bundles on the bounded symmetric domain $\Omega$. The homogeneous bundles can be obtained by holomorphic induction from representations of a certain parabolic Lie algebra on finite dimensional inner product spaces. The representations, and the induced bundles, have composition series with irreducible factors. In joint work with A. Kor\'{a}nyi 
\cite{KM0, KM}, our first main result is the construction of an explicit differential operator intertwining the bundle with the direct sum of its factors. Next, we study Hilbert spaces of sections of these bundles. We use this to get, in particular, a full description and a similarity theorem for homogeneous $n$-tuples of operators in the Cowen-Douglas class of the Euclidean unit ball in $\mathbb C^m$. A different approach is in \cite{MU}.
The initial study of these questions restricted to the case of homogeneous holomorphic line bundles is in \cite{BM}.

All these examples involve transitive action of the group $G$ on a homogeneous $G$- space, namely, $\Omega= G/H$ for some closed subgroup $H \subseteq G$. The study of homogeneity under the action of a group not necessarily acting transitively is more complicated. A small class of operators homogeneous under the action of the unitary group has been identified in \cite {GKMP}.

\begin{problem} Let $\mathbb{B}$ be the Euclidean ball in $\mathbb{C}^d$ and $\boldsymbol{T}$ be a commuting $d$- tuple of operators in the Cowen-Douglas class $\mathrm{B}_k(\mathbb{B})$.  The unitary group $\mathcal{U}(d)$ acts on any $d$- tuple of commuting operators as follows
\[u \cdot \boldsymbol{T}:=\big (\left(u_{11} T_1+ \cdots + u_{1d} T_d), \ldots, (u_{d1} T_1+ \cdots + u_{dd} T_d) \right), u \in \mathcal{U}(K).\]
We say that $\boldsymbol{T}$ is $\mathcal{U}(d)$- homogeneous if $u\cdot \boldsymbol{T}$ is unitarily equivalent to $\boldsymbol{T}$ for all $u\in \mathcal{U}(d)$.
Find all the $\mathcal{U}(d)$- homogeneous operators in  $\mathrm{B}_k(\mathbb{B})$.
\end{problem}

{\bf Algebraic geometry and operator theory.} The very successful Sz.-Nagy -- Foias model theory for contractions has two parts. First, one identifies a Hilbert module over the disc algebra which serves as an universal object and then obtain every contractive Hilbert module as a quotient of this universal object. Second, for any contractive module $(\mc{H}, T)$, Sz.-Nagy and Foias define a characteristic operator $\Theta_T$ and show that it is a complete unitary invariant for $(\mc{H}, T)$. The Nagy -- Foias theory builds on two fundamental results, the first is the von Neumann - Wold decomposition of an isometry and the other is the Beurling's characterization of the invariant subspaces of the shift operator.
Here is a special case: Any contractive Hilbert module $(\mc{H}, T)$ in the class $C_{\cdot 0}$ is a quotient of the Hardy module $H^2_{\mc E}(\mathbb D)$ by a submodule $\mc{S}$ which must be of the form $\Theta H_{\mc E}^2(\mathbb D)$ for some inner function $\Theta$ by Beurling's theorem. Moreover, $\Theta$ is a complete invariant for the module $(\mc{H}, T)$.

One proof of the Beurling's theorem
can be made up from the observation that all the submodules of the Hardy module are isomorphic. Interestingly, the quotient modules corresponding to these provide models for ``practically'' all the contractive modules. It is easy to produce examples in submodules of $H^2(\mb{D}^2)$ that are not of the form one would predict from the Beurling characterization. For example, take the submodule $\mc{S}_0$ consisting of functions vanishing at $(0,0)$. Not only that the submodule $\mc{S}_0$ is not isomorphic to $H^2(\mb{D}^2)$ as opposed to the case of one variable. On the other hand the submodules fail to detect the class of the quotient module. To see this, consider the submodules $\mc{S}_1$ and $\mc{S}_2$ of functions vanishing on the set $\{z_1=0\}$ and $\{z_2=0\}$ respectively. The quotient modules in this case are both isomorphic to the Hardy module $H^2(\mb{D})$.

One of the main difficulties here is that the dimension of the eigenspace of a submodule, in general (unless it is the completion of a principal ideal), is not necessarily constant. Thus these fail to define a holomorphic vector bundle. Indeed, for the submodule $\mc{S}_0$ of $H^2(\mb{D}^2)$, we have $\dim  \cap_{i=1,2} \ker M_i^* = 2$ while $\dim  \cap_{i=1,2} \ker (M_i - w_i)^* = 1$ if $(w_1, w_2)  \not = (0,0)$. Thus the map $(w_1, w_2) \to \cap_{i=1,2} \ker (M_i - w_i)^*$ does not define a holomorphic vector bundle on $\mb{D}^2$, however, it defines a coherent sheaf on $\mb{D}^2$ with a Hermitian structure inherited from that of $H^2(\mb{D}^2)$.

In general, determining the moduli space for the  isomorphism classes of sub-modules of a Hilbert module is a difficult problem. Thanks to Beurling's theorem, the moduli space for submodules of the Hardy module $H^{2}(\mathbb{D})$ is a singleton; that is, all submodules are isomorphic. However, a rigidity phenomenon occurs in $H^2(\mathbb D^n)$, namely,  no two sub-modules of $H^2(\mathbb D^n)$ are isomorphic barring a very few exceptions. This rigidity phenomenon is typical in multivariable settings.

Let $\Omega \subseteq \mathbb C^m$ be a bounded connected open set and $\mathcal H \subseteq \mathcal O(\Omega)$ be an analytic Hilbert module, that is, the Hilbert space $\mathcal H$ possesses a reproducing kernel $K$, the polynomial ring $\mathbb C[\underline{z}]\subseteq \mathcal H$ is dense and the point-wise multiplication induced by $p\in  \mathbb C[\underline{z}]$ is bounded on $\mathcal H$. We fix an ideal $\mathcal{I}\subseteq\mathbb{C}[\boldsymbol{z}]$ and let $[\mathcal{I}]$ denote its completion in $\mathcal{H}$. Let $X:[\mathcal I] \to \mathcal H$ be the inclusion map.
Thus we have a short exact sequence of Hilbert modules
\[
\begin{split}
0 \rightarrow [\mathcal I] \xrightarrow{X} {\mathcal H} \xrightarrow{\pi} \mathcal Q \rightarrow 0,
\end{split}
\]
where the module multiplication in the quotient $\mathcal Q:=[\mathcal I]^\perp$  is given by the formula  $m_p f = P_{[\mathcal I]^\perp} (p f),$ $p\in \mathbb C[\underline{z}],\,f\in \mathcal Q$. The analytic Hilbert module $\mathcal H$ defines a subsheaf
$\mathcal S^\mathcal H$ of the sheaf $\mathcal O(\Omega)$ of holomorphic functions defined on $\Omega$. For any open $U \subset \Omega$, it is obtained by setting
$$\mathcal S^\mathcal H(U) := \big \{\, \sum_{i=1}^n ({f_i|}_U) h_i : f_i \in \mathcal H, h_i \in \mathcal O(U), n\in\mathbb N\,\big  \}.$$
This is locally free and naturally gives rise to a holomorphic line bundle on $\Omega$. However, in general, the sheaf corresponding to the sub-module $[\mathcal I]$ is not locally free but only coherent.

Localising the modules $[\mathcal I]$ and $\mathcal H$ at $w\in \Omega$, we obtain the localization $X(w)$ of the module map $X$. The localizations are nothing but the quotient modules $[\mathcal I]/{[\mathcal I]_w}$ and $\mathcal H/{\mathcal H_w}$, where $[\mathcal I]_w$ and $\mathcal H_w$ are the maximal sub-modules of functions vanishing at $w$. These localizations define anti-holomorphic line bundles $E_{[\mathcal{I}]}$ and $E_{\mathcal{H}}$, respectively, on $\Omega\setminus V(\mathcal{I})$. However, there is a third line bundle, namely, ${\rm Hom}(E_\mathcal H, E_{[\mathcal I]})$ defined by the anti-holomorphic map $X(w)^*$. The curvature of a holomorphic line bundle $\mathcal L$ on $\Omega$, computed relative to a holomorphic frame $\gamma$ is given by the formula
$$\mathcal K_\mathcal L(z) =\sum_{i,j=1}^{m}\tfrac{\partial^2}{\partial z_i \partial \bar{z}_j}\log\|\gamma(z)\|^2 dz_i \wedge d\bar{z}_j.$$
It is a complete invariant for the line bundle $\mathcal L$. Now, consider the alternating sum
$$
\mathcal A_{[\mathcal I], \mathcal H}(w):=\mathcal K_X(w) - \mathcal K_{[\mathcal{I}]}(w) + \mathcal K_{\mathcal{H}}(w) = 0,\,\, w\in \Omega \setminus V_{[\mathcal I]},
$$
where $\mathcal K_X$, $\mathcal K_{[\mathcal{I}]}$ and $\mathcal K_{\mathcal{H}}$ denote the curvature $(1,1)$ form of the line bundles $E_X$, $E_{[\mathcal{I}]}$ and $E_{\mathcal{H}}$, respectively, on $\Omega \setminus V_{[\mathcal I]}$. Thus it is an invariant for the pair $([\mathcal I], \mathcal H)$.

\begin{problem}
The eigenspace distribution of an analytic Hilbert module $(\mc{H}, \bs{T})$ defines a coherent sheaf $\mathcal S^{\mathcal H}$ over $\Omega$ possessing a Hermitian structure. Find tractable unitary invariants for analytic Hilbert modules $(\mc{H}, \bs{T})$  from the sheaf model.
\end{problem}
The work along these lines was initiated in \cite{DMV1, G BM, BMP, DuanGuo},  however,  see \cite{BMS, Up1, Up2, Up3} for the most recent work on this topic.

{\bf Flag Structure.} Fix a bounded planar domain $\Omega$. Let $E$ be a holomorphic Hermitian vector bundle of rank $n$ in $\Omega \times \mathcal H$. By a well-known theorem of Grauert, every holomorphic vector bundle over a plane domain is trivial. So, there exists a holomorphic frame  $\gamma_1, \ldots , \gamma_n: \Omega \to \mathcal H$. A Hermitian metric on $E$ relative to this frame is given by the formula
$G_{\gamma}(w)= \left( \langle \gamma_i(w), \gamma_j(w)\rangle \right)$. The curvature of the vector bundle $E$ is a complex $(1,1)$ form which is given by the formula  \index{second fundamental form}
$$ \mathcal{K}_E(w)= \sum_{i,j=1}^m \mathcal{K}_{i,j}(w)\,  dw_i\wedge d\bar{w}_j,$$
where
$$\mathcal{K}_{i,j}(w) = -\frac{\partial}{\partial \bar{w}_j}\big (G_{\gamma}^{-1}(w) (\frac{\partial}{\partial w_i} G_{\gamma}(w))\big ).$$
It clearly depends on the choice of the frame ${\gamma}$ except when $n=1$. A complete set of invariants is given in
\cite{CD}. However, these invariants are not easy to compute. So, finding a tractable set of invariants for a smaller class of vector bundles which is complete would be worthwhile. For instance, in the paper \cite{JJKM}, irreducible holomorphic Hermitian vector bundles,  possessing a flag structure have been isolated. For these, the curvature together with the second fundamental form (relative to the flag) is a complete set of invariants. As an application, at least for $n=2$, it is shown that the homogeneous holomorphic Hermitian vector bundles are in this class. A complete description of these is then given. This is very similar to the case of $n=1$ except that now the second fundamental form associated with the flag has to be also considered along with the curvature.
All the vector bundles in this class and the operators corresponding to them are irreducible. The flag structure they possess by definition  is rigid which  aids in the construction of a canonical model and in finding a complete set of unitary invariants. The study of commuting tuples of operators in the Cowen-Douglas class possessing a flag structure is under way.

The definition of the smaller  class $\mathcal FB_2(\Omega)$ of operators in $B_2(\Omega)$ below is from \cite{JJKM}. One may similarly define $\mathcal FB_n(\Omega)$, $n >1$,  and discuss the properties of this new class, see \cite{JJKM}.

\subsubsection*{Definition} We let $\mathcal FB_2(\Omega)$ denote the set of bounded linear operators $T$ for which we can find operators  $T_0,T_1$  in $\mathcal{B}_1(\Omega)$ and  a non-zero intertwiner
$S $ between $T_0$ and $T_1$, that is, $T_0S=S T_1$ so that
$T=\begin{pmatrix}
 T_0 & S \\
 0 & T_1 \\
\end{pmatrix}.$

An operator $T$ in  $B_2(\Omega)$ admits a decomposition of the form  $\Big (\begin{smallmatrix}
T_0 & S \\
0 & T_1 \\
\end{smallmatrix}\Big )$
for some pair of operators $T_0$ and $T_1$ in ${B}_1(\Omega)$
(cf.  \cite[Theorem 1.49]{jw}) .
Conversely, an operator $T,$ which admits a decomposition of this form for some choice of $T_0, T_1$ in $B_1(\Omega)$ can be shown  to be in $B_2(\Omega).$
In defining the new class $\mathcal FB_2(\Omega)$, we are merely imposing one additional condition, namely that
$T_0S=ST_1$.

An operator $T$ is in the class $\mathcal FB_2(\Omega)$ if and only if there exist a frame
$\{\gamma_0,\gamma_1\}$ of the vector bundle $E_{T}$ such
that $\gamma_0(w)$ and $t_1(w):=\tfrac{\partial}{\partial w}\gamma_0(w)-\gamma_1(w)$ are orthogonal for all $w$ in $\Omega$.  This is also equivalent to the existence of a frame ${\gamma_{0}, \gamma_{1}}$ for the vector bundle $E_{T}$ satisfying $\frac{\partial}{\partial w}|\gamma_{0}(w)|^{2} = \langle\gamma_{1}(w), \gamma_{0}(w)\rangle$ for all $w \in \Omega$.

\subsubsection*{\rm Theorem  \cite[Theorem 2.10]{JJKM}}\label{maint}
Let
$T=\begin{pmatrix}T_0 & S \\
0 & T_1 \\
\end{pmatrix},\,\, \tilde{T}=\begin{pmatrix}\tilde{T}_0 & \tilde S \\
0 & \tilde{T}_1 \\
\end{pmatrix}$
be two operators in $\mathcal FB_2(\Omega)$. Also let $t_1$ and
$\tilde{t}_1$ be non-zero sections of the holomorphic Hermitian vector
bundles $E_{T_1}$ and $E_{\tilde{T}_1}$ respectively. The operators $T$ and
$\tilde{T}$ are equivalent if and only if
$\mathcal{K}_{T_0}=\mathcal{K}_{\tilde{T}_0}$ (or
$\mathcal{K}_{T_1}=\mathcal{K}_{\tilde{T}_1}$) and
$$\frac{\|S(t_1)\|^2}{\|t_1\|^2}= \frac{\|\tilde
S(\tilde{t}_1)\|^2}{\|\tilde{t}_1\|^2}.$$
Cowen and Douglas point out in \cite{CD} that an operator in $B_1(\Omega)$ must be irreducible.  However, determining which operators in $B_n(\Omega)$ are irreducible is a formidable task. Operators in $\mathcal{F}B_{2}(\Omega)$ are always irreducible. Indeed, if we assume $S$ is invertible, then $T$ is strongly irreducible.

Recall that an operator $T$  in the Cowen-Douglas class $B_n(\Omega)$,  up to unitary equivalence, is the  adjoint of the multiplication operator {$M$} on a Hilbert space {$\mathcal H$} consisting of holomorphic functions on $\Omega^*:=\{\bar{w}: w \in \Omega\}$ possessing a reproducing kernel {$K$.}  What about operators in $\mathcal FB_n(\Omega)$? For $n=2,$ a model for these operators is described below.

For an operator $T\in \mathcal FB_2(\Omega)$, there exists a  holomorphic frame $\gamma=(\gamma_0,\gamma_1)$ with the property $\gamma_1(w):=\tfrac{\partial}{\partial w}\gamma_0(w)-t_1(w)$ and that $t_1(w)$ is orthogonal to $\gamma_0(w)$, $w\in\Omega$, for some holomorphic map $t_1:\Omega \to \mathcal H$. In what follows, we fix a holomorphic frame with this property. Then the operator $T$ is unitarily equivalent to the adjoint of the multiplication operator $M$ on a Hilbert space
$\mathcal{H}_{\Gamma} \subseteq {\rm Hol}(\Omega^*, \mathbb C^2)$ possessing a reproducing kernel $K_{\Gamma}:\Omega^* \times \Omega^* \to \mathcal M_2(\mathbb  C)$. The details are in \cite{JJKM}. It is easy to write down the kernel $K_\Gamma$ explictly: For $z,w\in\Omega^*$, we have
\begin{eqnarray*}
K_{\Gamma}(z,w)&=&
\begin{pmatrix}
  \langle \gamma_0(\bar w),\gamma_0(\bar z)\rangle &
  \langle \gamma_1(\bar w),\gamma_0(\bar z)\rangle \\
   \langle \gamma_0(\bar w),\gamma_1(\bar z)\rangle &
   \langle \gamma_1(\bar w),\gamma_1(\bar z)\rangle \\
\end{pmatrix}\\
&=& \begin{pmatrix}
  \langle \gamma_0(\bar w),\gamma_0(\bar z)\rangle &
  \frac{\partial}{\partial \bar w}\langle \gamma_0(\bar w),\gamma_0(\bar z)\rangle \\
  \frac{\partial}{\partial z} \langle \gamma_0(\bar w),\gamma_0(\bar z)\rangle
  &
 \frac{\partial^2}{\partial z\partial \bar w} \langle \gamma_0(\bar w),
 \gamma_0(\bar z)\rangle+\langle t_1(\bar w),t_1(\bar z)\rangle \\
\end{pmatrix}, w\in\Omega.
\end{eqnarray*}
Setting $K_0(z,w)=\langle \gamma_0(\bar w),\gamma_0(\bar z)\rangle$
and $K_1(z,w)= \langle t_1(\bar w),t_1(\bar z)\rangle$, we see that the
reproducing kernel $K_{\Gamma}$ has the form:
\[ K_{\Gamma}(z,w)= 
\begin{pmatrix}
  {K_0}(z,w) & \frac{\partial}{\partial \bar w}{K_0}(z,w) \\
  \frac{\partial}{\partial z}{K_0}(z,w) & \frac{\partial^2}{\partial z\partial \bar w}{K_0}
  (z,w)+{K_1}(z,w) \\
\end{pmatrix}. \]

\begin{problem} Give an intrinsic definition of operators in the class $\mathcal{F}B_n(\Omega)$. Find a smaller class of commuting $d$- tuples of operators in $\mathrm B_n(\Omega)$ for a bounded domain $\Omega \subset \mathbb{C}^d$ possessing a tractable set of unitary invariants similar to the case of $n=1$.
\end{problem}
\subsubsection*{Final Remarks} This note does not address several important and interesting themes, such as:
\begin{enumerate}
\item Curvature inequalities (cf. \cite{BKM,GM1,Ramiz,GMRR,WZ,U,KT}),

\item Similarity problems within the Cowen-Douglas class (cf. \cite{JKSX,JJ,JJK,HJJX,KT}), and

\item Refinements of Cowen-Douglas operators exhibiting a flag structure (cf. \cite{XJ,YJ,JJM}).
\end{enumerate}

The references provided are intended to be representative, rather than exhaustive.
That said, I hope this brief account can help future work in these areas.

\section{Fourier theory, angle, Pimsner-Popa bases and regularity for inclusions of $C^*$-algebras}

\noindent \textit{By Keshab Chandra Bakshi [keshab@iitk.ac.in] from IIT Kanpur; and Ved Prakash Gupta [vedgupta@jnu.ac.in, vedgupta@mail.jnu.ac.in] from JNU, New Delhi.}

\bigskip

As is evident in  various categories, it is a widely accepted line of thought  in mathematics that the study of subobjects of a given object (from qualitative as well as quantitative perspectives) is one of the most effective  methodologies employed to obtain a good understanding of the structure and properties of the ambient object. {In fact, this approach  has led to some remarkable classification results in the categories of groups,  von Neumann algebras and $C^*$-algebras.}

Over the last five decades or so, various tools and techniques have been developed by various authors to study the relative positions between subalgebras of operator algebras.

This report focuses primarily on our recent research related to (i) developing the Fourier theory for inclusions of $C^*$-algebras; (ii) a new tool of angle between compatible intermediate operator subalgebras of finite-index inclusions of operator algebras; (iii) some questions of Jones and Popa regarding Pimsner-Popa bases (a heavily used tool  to study inclusions of operator algebras); and, (iv) obtaining a better understanding of the so-called regular inclusions of operator algebras, which depended heavily on the progress that we made under (iii).

\smallskip

\noindent {\bf $\bullet$ Angle and Fourier theory for inclusions of $C^*$-algebras }

Very recently, (along with some other authors)  we introduced the notions of interior and exterior angles between compatible intermediate subalgebras of finite-index inclusions of certain operator algebras (see \cite{BDLR, BG3}). Here is a very brief report on this development.

Vaughan Jones (\cite{Jones}) introduced the notions of index and  basic construction for any unital inclusion $N \subset M$ of  $II_1$-factors. Later, Bisch (\cite{Bisch}) showed that there exists a dictionary between intermediate subfactors of a finite-index inclusion $N \subset M$ and the so-called biprojections in the relative commutant of $N$ in the basic construction $ M_1$. Using this correspondence, the notions of interior and exterior angles between any two intermediate subfactors of a finite-index inclusion $N \subset M$ was introduced in \cite{BDLR}. As an application, when a finite-index inclusion $N \subset M$ is irreducible (i.e., $N'\cap M = \mathbb{C}$), employing some aspect of the (well-known) Fourier theory on the relative commutants of the subfactor, a question of Longo (\cite{Longo}) was answered positively by improving his bound for the cardinality of the lattice of intermediate subfactors of $N \subset M$. Later, some explicit calculations of angle between intermediate subfactors of group-subgroup subfactors were obtained in \cite{BG1}.

Motivated by \cite{BDLR}, using Watatani's theories of index and basic construction (from \cite{Wat}), the notions of interior and exterior angles between the so-called compatible intermediate $C^*$-subalgebras (as in \cite{IW}) of a finite-index inclusion $B \subset A$ of unital $C^*$-algebras was introduced recently in \cite{BG3} (and some calculations were done in \cite{GS}). (Following Watatani, the underlying theory for this line of research has been the notion of Hilbert $C^*$-modules.)  Interestingly, when $B \subset A$ are both simple, then there is a unique conditional expectation $E_0: A \to B$ with minimal Watatani index. Then, like the Jones' basics construction tower, one obtains an increasing tower of simple unital $C^*$-algebras $\{A_k: k \geq 1\}$.  Motivated by the theory of subfactors, we developed a Fourier theory on the relative commutants $\{ B'\cap A_k, A'\cap A_k : k \geq 0\}$ in terms of shift and rotation operators, which was enhanced further in \cite{BGS, BGP}. As an application, in \cite{BG3} itself, we could provide a bound for the cardinality of the
lattice of intermediate $C^*$-subalgebras of a finite-index irreducible inclusion of simple unital $C^*$-algebras (whose finiteness was already proved by Ino and Watatani (\cite{IW}); and, we
could improve Longo's bound for the lattice of intermediate subfactors of finite-index irreducible inclusions of type $III$ factors as
well.

\smallskip

\noindent {\bf $\bullet$ Pimsner-Popa bases and regularity for subfactors and inclusions of simple unital $C^*$-algebras}

Over the last four decades or so, the notion of the so-called Pimsner-Popa basis (\cite{PiPo}) for finite-index inclusions of $II_1$-factors has played an indispensable role in the development of the theory of subfactors (via Fourier theory, Popa's $\lambda$-lattices and Jones' planar algebras, to name a few). Jones was particularly impressed by their utility and went on to ask (on
more than one occasion) whether every so-called extremal inclusion of $II_1$-factors always admits a two-sided Pimsner-Popa basis or not. This question remains open till date. However, in \cite{BG2}, we provided a partial answer to Jones' question by proving that a finite-index regular inclusion of $II_1$-factors always admits a two-sided Pimsner-Popa basis and has integer index. (An inclusion $P
\subset Q$ of unital $C^*$ or $W^*$-algebras is said to be regular if the unitary normalizers of $P$ in $Q$ generate $Q$ in the $C^*$ or $W^*$ sense.) This was achieved by an appropriate application of the notion of path algebras developed independently by Sunder and Ocneanu.

Recently, Popa (\cite{Popa}) asked whether every integer-index irreducible inclusion of $II_1$-factors admits an orthonormal (Pimnser-Popa) basis or not. In \cite{BG4}, we gave a positive answer to Popa's question for finite-index regular subfactors with simple or commutative relative commutants. These results then allowed us to deduce (from \cite{NV}) that such a subfactor is isomorphic to a crossed-product subfactor
with respect to a minimal action of a so-called biconnected weak Kac algebra on a $II_1$-factor.

Based on the progress made in \cite{BG4} and some techniques from the world of Quantum Information theory, Crann et al (\cite{CKP}) recently showed that every finite-index regular subfactor admits an orthonormal (Pimsner-Popa) basis and thereby answering a conjecture of ours from \cite{BG4}. As a result, a proof from \cite{BG4} implies that every finite-index regular subfactor has depth at most 2. In the finite-dimensional setting, for appropriate inclusions of finite-dimensional $C^*$-algebras, Popa's question has been further investigated in \cite{BB} (once again) by employing some techniques from Quantum Information Theory.

Over the last 3 decades or so operator algebraists have been actively trying to find similarities between the theory of subfactors and the theory of inclusions of simple $C^*$-algebras. From the works of various authors, it is well-known that every irreducible finite-index subfactor is isomorphic to a crossed-product subfactor with respect to an outer action of a finite group. Motivated by this characterization, quite satisfyingly, we have recently been able to show (in \cite{BG5}) that every finite-index regular inclusion of simple unital $C^*$-algebras is isomorphic to a reduced cocycle crossed-product inclusion with respect to an ``outer'' action of a finite group. Moreover, we were also able to prove that any such inclusion of simple unital $C^*$-algebras has integer (Watatani) index, has depth at most 2 and the minimal conditional expectation of such an inclusion admits a two-sided Pimsner-Popa basis.

\section{Contributions in Banach space theory}

\noindent \textit{By Pradipta Bandyopadhyay [pradipta@isical.ac.in] from ISI Kolkata.}

\bigskip

\noindent \textsf{On $L_1$-predual spaces:} A Banach space $X$ is called an $L_1$-predual space or a Lindenstrauss space if $X^* = L_1(\mu)$, for some measure $\mu$.

In \cite{DKS}, motivated by Bratteli diagrams of Approximately Finite Dimensional (AF) $C^*$-algebras, the authors consider diagrammatic representations of separable $L_1$-predual spaces and show that, in analogy to a result in AF $C^*$-algebra theory, in such spaces, every directed sub-diagram represents an $M$-ideal. The converse, namely, the question whether given an $M$-ideal in a separable $L_1$-predual space $X$, there exists a diagrammatic representation of $X$ such that the $M$-ideal is given by a directed sub-diagram, remains open in general. We refer to this as the ``main problem''. In \cite{DKS}, the authors give affirmative answers to this in some special cases.

Continuing this study, we have shown in \cite{BS1} that for a compact metric space $K$, in $C(K)$-space, our ``main problem'' has an affirmative answer. Similar results are true also for the spaces $C_0(T)$, $\textbf{c}$ and $\textbf{c}_0$.

In \cite{BDS2}, we present characterizations of polyhedrality $(I)$, $(II)$ and $(III)$ in a separable $L_1$-predual space in terms of its representing matrix. We also show that in a polyhedral $(IV)$ predual of $L_1$, our ``main problem'' has an affirmative answer.

In \cite{BDS1}, we study non-separable Gurariy spaces and their almost isometric ideals. The main result is that a (non-separable) Banach space is a Gurariy space if and only if every separable \emph{almost isometric ideal} (a.i.-ideal) in $X$ is isometric to the separable Gurariy space $\mathbb{G}$. We also obtain a similar characterization of $L_1$-predual spaces in term of \emph{ideals}. Along the way, we show that the family of ideals/a.i.-ideals in a Banach space is closed under increasing limits. And hence, the family of all separable ideals/a.i.-ideals in a Banach space is a skeleton.

Continuing in \cite{BS2}, we show that if a Banach space $X$ is almost isometric to a Gurariy space, then it is also a Gurariy space. In the separable case, this answers a question posed by Rao \cite[Remark 6]{R1} in the affirmative. We also prove an analogue of this result for the noncommutative Gurariy space. We also show that a Banach space $X$ is a Gurariy space if and only if it is an a.i.-ideal in every superspace in which it embeds as a hyperplane, thereby answering another question posed by Rao \cite[Question 17]{R} in the affirmative.

\smallskip

\noindent \textsf{On Mazur Intersection Property and its variants:} The Mazur Intersection Property (MIP)---every closed bounded convex set is the intersection of closed balls containing it---is an extremely well studied property in Banach space theory.

A complete characterisation of MIP was obtained in Giles, Gregory and Sims \cite{GGS}, most well-known criterion stating that the w*-denting points of $B(X^*)$ are norm dense in $S(X^*)$. The paper also considered w*-MIP in dual spaces.

A much less studied uniform version of MIP (UMIP) was introduced by Whitfield and Zizler \cite{WZ}. Characterisations similar to \cite{GGS} were also obtained, but an analogue of the w*-denting point criterion was missing, which perhaps is a reason for its being less pursued. A long standing open question regarding UMIP is whether UMIP implies that $X^*$ has a uniformly convex renormimg.

Chen and Lin \cite{CL} introduced the notion of w*-semidenting points and showed that a Banach space $X$ has MIP if and only if every $f \in S(X^*)$ is a w*-semidenting point of $B(X^*)$.

In \cite{BGG}, we show that a Banach space $X$ has UMIP if and only if every $f \in S(X^*)$ is a uniformly w*-semidenting point of $B(X^*)$, thus filling a long felt gap. In the process, we obtain simpler proofs of some characterisations in \cite{WZ}.

In \cite{Go}, we introduce two moduli of w*-semidenting points and characterise MIP and UMIP in terms of these moduli. One of the moduli was motivated by the modulus of denting points discussed by Dutta and Lin \cite{DL}. We introduce a property called Hyperplane UMIP (H-UMIP), which is slightly stronger than UMIP, and show that it characterises uniform convexity of $X^*$. This turns out to be also equivalent to the notion of uniform denting points of \cite{DL}. In this paper, we also analyse conditions, in terms of these moduli, for the stability of UMIP under $\ell_p$-sums, $1 <p <\infty$. In particular, we show that $\ell_p (X)$ has UMIP if and only if $X$ has UMIP.

Given a family $\mathcal{C}$ of closed bounded convex sets in a Banach space $X$, we say that $X$ has $\mathcal{C}$-MIP if every $C \in \mathcal{C}$ is the intersection of the closed balls containing it. This has also been studied in the literature for specific families and also in some generality.
Instead of considering just a family $\mathcal{C}$, in \cite{pBG1}, we discuss representation of sets of the form $\overline{A+\lambda B(X)}$, for $A \in \mathcal{C}$ and $\lambda \geq 0$, as intersection of closed balls. It is easy to see that this is stronger than $\mathcal{C}$-MIP. We call this strong $\mathcal{C}$-MIP and show that it is a more satisfactory generalisation of MIP inasmuch as one can obtain complete analogues of various characterisations of MIP in \cite{GGS}. We also define and characterise uniform versions of (strong) $\mathcal{C}$-MIP. Even in this case, strong $\mathcal{C}$-UMIP appears to have richer characterisations than $\mathcal{C}$-UMIP.

We say that a property $P$ of a Banach space is residual if it satisfies the following: If $X$ has one norm with property $P$, then ``almost every'' (residual with respect to some suitable metric) equivalent norm on $X$ also has property $P$. In \cite{pBG2}, we prove the residuality of UMIP, uniformly smooth and asymptotically uniformly smooth norms. Some of these are already known to be ball separation properties. For others, we first obtain a ball separation characterisation and then use them to obtain the residuality results.

It is well-known that Fr\'echet smoothness implies MIP. And since Fr\'echet smoothness is hereditary, it actually implies hereditary MIP, that is, all subspaces of $X$ also have MIP. How about the converse? Does hereditary MIP imply Fr\'echet smoothness? It is easy to see that hereditary MIP implies smoothness. In \cite{pBG3}, we obtain some sufficient conditions for hereditary MIP in terms of the set $N$ of non-Fr\'echet smooth  points of $S(X)$. In particular, we show that if $X$ is smooth and $N$ is finite, then $X$ has hereditary MIP.

Borwein and Fabian in \cite{BF} proved an infinite dimensional WCG Asplund space $X$ admits an equivalent smooth norm such that $N = \{\pm x_0\}$ for some $x_0 \in S(X)$. Thus, the Borwein-Fabian norm has hereditary MIP, but is not Fr\'echet smooth.

\section{Quantum dynamical maps and quantum Gaussian states}

\noindent \textit{By B V Rajarama Bhat [bhat@isibang.ac.in, bvrajaramabhat@gmail.com] from ISI, Bangalore.}
\bigskip

Here is  a  brief account of developments in recent years, through contributions by Indian mathematicians, in some mathematical areas that I am familiar with.  This text had to be prepared in a rather short period of time. Unfortunately this means that  majority of the papers cited here are from my collaborators and myself.

\smallskip

\noindent \textsf{Quantum Gaussian states:} Quantum Gaussian states on Boson Fock spaces are quantum versions of Gaussian distributions of classical probability.
The concept has been beautifully explained in the essay {\em ``What is a Gaussian state?''} by K. R. Parthasarathy \cite{Pa}.  In his late years he worked vigorously on this topic and discovered various features of convexity and symmetry properties of  finite mode quantum Gaussian states. He posed the problem of extending the theory to the infinite mode. This was carried out in \cite{R BJS}. It is seen that every quantum Gaussian state comes with
a `covariance matrix', which  is now a symmetric and invertible real linear
operator. A complete characterization of covariance operators of quantum
Gaussian states in infinite mode has been obtained  and it involves some delicate trace conditions which are not seen in finite mode.
A unitary operator in the Boson Fock space is called a {\em Gaussian symmetry\/}  if
it preserves Gaussianity of states under conjugation. It is shown
that the symplectic spectrum is a complete invariant for Gaussian states under conjugation by Gaussian symmetries.

The equivalence of separability and complete extendability for bipartite quantum Gaussian states is seen in \cite{BPS}. Some of the questions raised by K R Parthasarathy have been answered recently in \cite{DJS}. In \cite{TP} we can find a novel new parametrization for quantum Gaussian states. It is very useful in some of the computations involving these states \cite{Pa2} and perhaps its full potential hasn't yet been exploited.

As mentioned above, while dealing with quantum gaussian states the symplectic spectrum plays an important role. This has triggered interest in studying approximation and perturbation properties of   symplectic eigenvalues (\cite{BTS}, \cite{LAST}).

\smallskip

\noindent \textsf{$C^*$-convexity of completely positive maps:} There have been several approaches to generalize the notion of convexity to the quantum context. One such  idea is to replace positive scalars in the interval $[0, 1]$ as coefficients for convexity by positive, contractive and invertible elements in a $C^*$-
algebra.  A study of the $C^*$-convexity structure of normalized positive operator
valued measures (POVMs) on measurable spaces has been carried out in \cite{BBK}. One of the surprising
results is that $C^*$-extreme points of normalized POVMs on
countable spaces (in particular for finite sets) are always spectral measures (normalized
projection valued measures). This does not hold if one considers usual convexity!
A Krein–Milman type theorem for POVMs has also been proved.

The set of all unital completely positive (UCP) maps on a unital $C^*$-algebra taking values in the algebra of all bounded operators on a Hilbert space is known as {\em generalized state space} and its convexity structure is of considerable interest. It is also a $C^*$-convex set.
The article \cite{BM1} links  $C^*$-extremity of this set with factorization property of certain algebras. Here also a Krein-Milman type result can be proved. This work brings forth an important connection between  $C^*$-convexity   theory and some classical concepts from operator algebra theory,  like nest algebras (infinite dimensional analogues of upper triangular matrices) and logmodular algebras.  A  not necessarily $*$-closed subalgebra $\mathcal{A} $  of a $C^*$-algebra $\mathcal{M} $ is said to be {\em logmodular}  (resp. has factorization) if the set
$\{a^*a: a \in \mathcal{M}~\mbox{is invertible with}~  a, a^{-1}\in \mathcal{A} \} $ is dense in (resp. equal to) the set of all positive
and invertible elements of $\mathcal{M}.$ One major result in \cite{BM2} is that the lattice of projections
in a (separable) von Neumann algebra $\mathcal{M}$ whose ranges are invariant under a logmodular
algebra in $\mathcal{M}$, is a commutative subspace lattice. Further, if $\mathcal{M}$ is a factor then this lattice
is a nest. As a special case, it follows that all reflexive  logmodular subalgebras of type I factors are nest algebras, and this settles a question which has been open for over a decade.  As a special case a  complete characterization of logmodular subalgebras in finite-dimensional von Neumann algebras has been obtained.

Recently the concept of $C^*$-convexity in the context of entanglement breaking (EB) maps has been studied in \cite{BDNS}. This has implications in quantum information theory.

\smallskip

\noindent \textsf{Peripheral Poisson boundary:} Unital completely positive (UCP) maps are quantum analogues of stochastic maps of classical probability. Suppose $\mathcal{A}$ is a von Neumann algebra and $\tau: \mathcal{A}\to \mathcal{A}$ is a normal UCP map. The space of fixed points of $\tau $ is a vector space but in general it is not a subalgebra of $\mathcal{A}$. However, it can be given a $*$-algebra structure with a modified product called Choi-Effros product. Taking motivation from classical probability, M. Izumi called this as the noncommutative Poisson boundary of the UCP map. The  noncommutative Poisson boundary has been computed for some specific Markov maps on the full Fock space in \cite{BBDR} and it produces some interesting von Neumann algebras.

   It was noted in \cite{BKT1} that Choi-Effros product can be extended to the closed linear span of the space of peripheral eigenvectors of $\tau$:
$$E_{\lambda }(\tau )=\{ x\in \mathcal{A}: \tau(x)=\lambda x\}, \lambda \in \mathbb{C}, |\lambda |=1.$$
 This way we get a $C^*$-algebra  called  the {\em peripheral Poisson boundary} of $\tau .$ Surprisingly it may not be possible to get a von Neumann algebra. The construction used {\em dilation theory} of completely positive maps developed earlier in \cite{Bh} and related papers.
This concept is very important as the asympotic behaviour of quantum dynamical semigroups are determined by peripheral eigenspaces.

\smallskip

\noindent \textsf{Iterative roots of functions:} Consider a self-map $f$ on a non-empty set. A
self-map $g$ is said to be an iterative $n$ -th root of $f$ if $f=g^n$, where the power is taken with respect to composition. A powerful and widely applicable new method has been developed in (\cite{R BG1}-\cite{R BG3}) to prove non-existence of such  roots under very general situations. For instance if there is distinguished point $x_0$, such that $f(x_0)\neq x_0$, $f^{-2}(\{x_0\})$ is infinite and $f^{-1}(\{x\})$ is finite for $x\neq x_0$, then $f$ has no iterative roots of any order $n\geq 2$.
This method can be   used to prove that continuous self-maps with no iterative $n$-th roots (for fixed $n>1$)  are dense in the space of all continuous self-maps for various topological spaces such as the unit cube in $\mathbb{R}^m$.  Existence and non-existence of iterative roots is of interest even in the theory of stochastic maps and quantum dynamical semigroups
\cite{R BHMV}.

\section{On Carath\'eodory approximation}

\noindent \textit{By Tirthankar Bhattacharyya [tirtha@iisc.ac.in] from IISc; and Poornendu Kumar [poornendukumar@gmail.com] from University of Manitoba, Canada.}

\bigskip

Approximation by polynomials or rational functions for holomorphic maps on the open unit disc $\mathbb{D}$ does not, in general, ensure any nice boundary behavior. Motivated by this limitation, Carath\'eodory, in 1926, proved the following remarkable theorem:
	
\medskip
	
\textbf{Theorem (Carath\'eodory, \cite{Caratheodory}.)} \emph{Any holomorphic function $\varphi : \mathbb{D} \to \overline{\mathbb{D}}$ can be approximated uniformly on compact subsets of $\mathbb{D}$ by the finite Blaschke products.}
	
\medskip
	
This theorem is a powerful tool in both holomorphic function theory and operator theory. For example, it can be used to prove the von Neumann inequality \cite{Drury}, derive the Herglotz integral representation \cite{Bhattacharyya-Bhowmik-Kumar}, and establish the Halmos conjecture—or more generally, the Berger--Stampfli theorem; see, for instance, \cite[Section 12]{Mashreghi-Rasnford}.
	
When seeking generalizations of this theorem, one must note that Carath\'eodory’s original proof relies heavily on the power series expansion of holomorphic functions on $\mathbb{D}$, which does not extend naturally to matrix-valued or operator-valued settings, nor to more complicated domains such as multiply connected ones. It is not difficult to observe that the finite Blaschke products are precisely the rational inner functions. Here, an \emph{inner function} refers to a bounded holomorphic map on $\mathbb{D}$ whose boundary values have modulus one almost everywhere with respect to the arc length measure on the unit circle $\mathbb{T}$. This notion extends naturally to the matrix-valued or operator-valued setting, where the modulus-one condition is replaced by the requirement that the boundary values are isometries almost everywhere on $\mathbb{T}$.
	
There is a proof of Carath\'eodory’s approximation theorem via the Pick–Nevanlinna interpolation theorem. This approach has the advantage of extending naturally to the matrix-valued case; see the discussion in the introduction of \cite{Alpay-Bhattacharyya-Jindal-Kumar}. However, this depends on the existence of inner solutions to the Pick-Nevanlinna interpolation problem, and such results are not known on multiply connected domains. Thankfully, one can deduce the Carath\'eodory approximation result on multiply connected domains from Grunsky’s work \cite{Grunsky}, and a different version exists due to Fisher \cite{Fisher2}. However, these approaches do not extend to matrix-valued settings. Hence, it requires new techniques. Recently, using the Herglotz integral representation and extreme point theory in the certain class, see \cite{Bhattacharyya-Bhowmik-Kumar}, this approximation theorem has been extended to the matrix-valued case. Moreover, taking into account the discussion at the beginning, we have that the Carath\'eodory approximation theorem and the Herglotz integral representation are essentially equivalent.

Furthermore, for a nice class of domains, it has been shown that once a matrix-valued Carath\'eodory approximation result is established \cite[Theorem 4.4]{Bhattacharyya-Bhowmik-Kumar}, it can be lifted to the operator-valued setting. Therefore, it suffices to focus on the matrix-valued version of the approximation theorem. As a result, operator-valued Carath\'eodory approximation theorems have now been obtained for multiply connected domains.
	
Moving to the several-variable setting, the scalar-valued Carath\'eodory approximation theorem has been established on the open unit ball independently by  Aleksandrov \cite{Aleksandrov} and Rudin \cite{Rudin2}. Moreover, Rudin extended this result to the polydisc \cite{Rudin}. More recently, the theorem has been proved for tube domains \cite{Wang-Zhang} and quotient domains related to the polydisc \cite{Bhowmik-Kumar}, domains arising from the classification of pseudoreflection groups.
	
However, existing techniques do not directly extend these results to matrix-valued settings in several variables. That said, through the use of the Nevanlinna-Pick interpolation theorem, it is possible to prove a matrix-valued analogue in the bidisc $\mathbb{D}^2$. In recent work \cite{Alpay-Bhattacharyya-Jindal-Kumar}, this theorem has been proved using operator theory, more specifically, dilation theory and realization formulas. Moreover, several applications of matrix-valued Carath\'eodory approximation results have been provided. To mention a few, Carath\'eodory’s theorem for matrix-valued functions and Potapov’s elegant characterization of matrix-valued rational inner functions on $\mathbb{D}$ \cite{Potapov} have led to a generalization of Fisher’s theorem \cite{Fisher}---namely, that a matrix-valued function which is holomorphic in $\mathbb{D}$ and continuous on $\overline{\mathbb{D}}$ can be uniformly approximated on $\mathbb{D}$ by convex combinations of matrix-valued rational inner functions---as well as to a version of Carath\'eodory’s theorem for functions taking values in the symmetrized bidisc or the tetrablock.

To include one proof in this note, we have chosen the proof of Carath\'eodory's theorem using dilation theory on the disc. The same proof works for the bidisc.  Given a holomorphic function $f : \mathbb{D} \to M_d(\mathbb{C})$ with $\|f(z)\| \leq 1$ for all $z\in\mathbb{D}$, one may first approximate it by rational functions. If $f$ is rational, then it admits a {\em finite contractive realization}, i.e.,
$
f(z) = A + z B (I - z D)^{-1} C,
$
where the associated colligation contraction matrix
$
T = \begin{bmatrix}
	A & B \\
		C & D
\end{bmatrix} \quad \text{acts on } \mathbb{C}^d \oplus \mathbb{C}^m
$
for some $m$. When $T$ is unitary, the resulting function $f$ is rational inner. This realization is a disc-version of the Kalman--Yakubovich--Popov lemma, see \cite{Dickinson-Delsarte-Genin-Kamp} for an indefinite version of it,  and its analogue in $\mathbb{D}^2$ has been recently proved by Knese \cite[Proposition 4.2]{Knese}.

Let $D_{T^{*}}$ and $D_{T}$ be the defect operators $(I- T T^{*})^{1/2}$ and $(I - T^{*} T)^{1/2}$ respectively. For each $N$, consider the space $\mathcal K_{N+1} = (\mathbb{C}^{d} \oplus \mathbb{C}^{m}) \oplus \dots \oplus (\mathbb{C}^{d} \oplus \mathbb{C}^{m}),$ where $(\mathbb{C}^{d} \oplus \mathbb{C}^{m}) $ occurs $(N+1)$ times and the {\em $N$-step unitary dilation of the contraction $T$} (first done by O.
	Nevanlinna \cite{Nevanlinna} and later greatly popularized by Levy and Shalit  \cite{Levy-Shalit}) given by
		$$ U_{N} :=  \begin{bmatrix}
			T & 0 &\dots & \dots & D_{T^*} \\
			D_T & 0 & \dots & \dots &  -T^*\\
			0 & I_{\mathcal{H}} & \dots  & \dots& 0 \\
			\vdots & \vdots & \ddots & \vdots & \vdots \\
			0 & 0 & \dots & I_{\mathcal{H}} & 0
		\end{bmatrix}  = \begin{bmatrix}
			A & B_{N} \\ C_{N} & D_{N}
		\end{bmatrix}$$
acting on the space $\mathcal K_{N+1}$ where $B_{N}, C_N$ and $D_N$ can be appropriately read off. Consider the rational inner function
$f_N(z) = A + z B_N (I - z D_N)^{-1} C_N$. A straightforward calculation shows that
	$
	f_N \to f \quad \text{if and only if} \quad A_N \to A \quad \text{and} \quad B_N D_N^k C_N \to B D^k C \quad \text{for each } k.
	$
The proof is complete because $ B_{N} D_{N}^{k}C_{N} = B D^{k} C $ for all $N \geq k+2$.
		
Knese \cite{Knese2} used the scalar-valued Carath\'eodory theorem on $\mathbb{D}^2$ to establish the necessity in the Nevanlinna–Pick interpolation problem on $\mathbb{D}^2$. We conclude with an open problem.
		
\quad
	
\textbf{Open Problem.} \emph{Does a Carathéodory approximation theorem hold for matrix-valued functions on the polydisc?}

\textsl{The first named author is supported by ANRF(JCB/2021/000041) and the second named author is supported by a PIMS post-doctoral fellowship.}

\section{Operator theory on quotient domains}

\noindent \textit{By Shibananda Biswas [shibananda@iiserkol.ac.in] and \\
Subrata Shyam Roy [ssroy@iiserkol.ac.in] from IISER Kolkata.}

Models for $n$-tuples of operators with the polydisc and the Euclidean ball in $\mathbb C^n$ as spectral sets have been studied in depth. A natural progression is to examine such models on inhomogeneous domains \cite[Corollary 4.4]{AY2}. The class of $\Gamma_n$-contractions provided the first instance in this pursuit. Agler and Young \cite{AY, AY1} found that many of the fundamental results (e.g substantial part of the Sz-Nagy-Foias model theory) in the theory of single contraction has close parallels for $\Gamma_2$-contractions - commuting pair of operators with closed symmetrized bidisc $\Gamma_2$ as a spectral set. Motivated by these developments, our attention turned to functional models on (open) symmetrized polydisc $\mathbb G_n$ \cite{BS}, beginning with the construction of a Schaffer dilations for $\Gamma_2$-contractions \cite{S BPS} and exploring possible generalization to $n$-tuples. A  Beurling-Lax-Halmos type theorem for joint invariant subspaces of a pure $\Gamma_n$-isometry is obtained by generalizing the same in \cite{S} for $n=2$ case. In fact on $\Gamma_2$,  a number of questions were answered satisfactorily, for example, description of distinguished variety \cite{S PS}, on realization, interpolation and extension of holomorphic functions \cite{BhattS, ALY, S BKS}. The generalization of several results in this direction beyond the case $n=2$ has been hindered by their dependence on dilation theory as positive results towards rational dilation of $\Gamma_n$-contraction, for $n>2$, remained elusive beyond restricted subclasses \cite{S AP}. A graph-theoretic model for a large class of $n$-tuple of operators that includes $\Gamma_n$-contractions as a special case was considered in \cite{GH}.

 On the other hand, in \cite{MSZ}, function theory on $\mathbb G_n$ was studied by introducing a Hilbert space method to compute the Bergman kernel of $\mathbb G_n$ which was earlier obtained in \cite{EZ} appealing to  Bell's formula \cite{B}. This work naturally led to the consideration of Hilbert function spaces on $\mathbb G_n$ and the study of operators acting upon them. In \cite{BS}, we constructed functional models for $\Gamma_n$-contractions on weighted Bergman spaces $\mathbb A^{(\lambda)}(\mathbb D^n)$, $\lambda>1$, and on the Hardy space $H^2(\mathbb D^n)$ as the limiting case $\lambda = 1$. These spaces decompose into orthogonal direct sum of submodules over the invariant ring $\mathbb C[\boldsymbol z]^{\mathfrak S_n}$, each of which corresponds to an irreducible representation of the
symmetric group $\mathfrak S_n$ and hence, indexed by partitions $\boldsymbol p$ of $n$ - denoted as $\boldsymbol p\vdash n$. The restrictions of $(M_{s_1},\ldots,M_{s_n})$ to these submodules are $\Gamma_n$-contractions. In \cite{BGMS}, this framework was generalized to Hilbert spaces with $\mathfrak S_n$-invariant reproducing kernels $K$, leading to reducing submodules $\mathbb P_{\boldsymbol p}\mathcal H$ with kernel
\(
  {\rm imm}^{(\boldsymbol p)}\big(\!\!(K(z_i,w_j))\!\!\big)_{i,j=1}^n,
\)
where ${\rm imm}^{(\boldsymbol p)}$ is the immanant corresponding to a partition $\boldsymbol p$ and  the map $\mathbb P_{\boldsymbol p}:\mathcal H \to \mathcal H$ is an orthogonal  projection given by the formula
$$
  \mathbb{P}_{\boldsymbol p} f  = \frac{\deg {\boldsymbol p}}{n!}\sum_{\sigma\in\mathfrak S_n}{\chi_{\boldsymbol p}(\sigma^{-1})} (f \circ \sigma^{-1}),\,\,f\in \mathcal H,
 $$
 $\chi_{\boldsymbol p}$ being the character of the representation corresponding to $\boldsymbol p\vdash n$. A finer decomposition into the reducing submodules $\mathbb P_{\boldsymbol p}^{ii}\mathcal H$, $1\leq i\leq \deg \boldsymbol p$, was also obtained. This naturally led to the study of unitary equivalence, similarity, and minimality of these reducing submodules as well as their realization on $\mathbb G_n$.

\emph{\textbf{On unitary equivalence.}} Each $\mathbb P_{\boldsymbol p}\mathcal H$ is a locally free Hilbert module of rank $(\deg \boldsymbol p)^2$. Thus, if $\deg \boldsymbol p \neq \deg \boldsymbol q$, then the  submodules $\mathbb P_{\boldsymbol p}\mathcal H$ and $\mathbb P_{\boldsymbol q}\mathcal H$ are not unitarily equivalent. When $\deg \boldsymbol p = \deg \boldsymbol q$, the problem remains open in general. However, if $\deg \boldsymbol p=\deg \boldsymbol q=1$ and $\mathcal H=\mathbb A^{(\lambda)}(\mathbb D^n)$ for  the reducing  submodules corresponding to  the trivial and sign representations of $\mathfrak S_n$ are not equivalent for $\lambda >0$ \cite{BGMS}. This follows from the identification of the submodule of $\mathbb A^{(\lambda)}(\mathbb D^n)$ corresponding to the sign representation of $\mathfrak S_n$ as the weighted Bergman module $\mathbb A^{(\lambda)}(\mathbb G_n),$ and by obtaining a formula for the reproducing kernel of $\mathbb A^{(\lambda)}(\mathbb G_n)$ in the coordinates of $\mathbb G_n$ invoking Giambelli's formula which is of independent interest.

\emph{\textbf{On similarity.}} Even when unitary equivalence fails, the possibility of similarity remains. If $\boldsymbol p = (n)$ and $\boldsymbol q = (1^n),$ then a candidate for implementing such an isomorphism is the operator of multiplication by Vandermonde $\Delta$,  where $\Delta(\boldsymbol z)=\prod_{1\leq i<j\leq n}(z_i-z_j)$, from $\mathbb P_{\boldsymbol p}\mathcal H$ to $\mathbb P_{\boldsymbol q}\mathcal H$. Multiplication by the Vandermonde is injective, therefore it remains only to check if it is surjective to make it invertible.  We rephrase the question of surjectivity in the realm of a division problem: If for $f \in \mathbb P_{\boldsymbol q}\mathcal H$, suppose there is a $g$ such that  $ \Delta g = f$, then does it follow that $g$ is in $\mathbb P_{\boldsymbol p}\mathcal H$? It would be really interesting to know, as the results that are known mostly in the case of  strictly pseudoconvex domains, extend to the case of the product domain $\mathbb D^n,$ even in the degree $1$ case (cf. \cite{PuS} and references therein).

\emph{\textbf{On minimal reducing subspaces.}} For each partition $\boldsymbol p \vdash n$ and each index $i$ with $1 \leq i \leq \deg \boldsymbol p$, the submodules $\mathbb P_{\boldsymbol p}^{ii}\mathcal H$ are reducing for the $n$-tuple $(M_{s_1},\ldots,M_{s_n})$. A natural question, therefore, is whether these submodules are minimal reducing. A more general perspective considers a finite complex reflection group $G$ with homogeneous generators $\boldsymbol\theta=(\theta_1,\ldots,\theta_n)$ of the invariant ring $\mathbb C[\boldsymbol z]^G.$ Within this framework, multiplication by $\boldsymbol \theta$ on a  Hilbert module $\mathcal H$ with $G$-invariant reproducing kernel over a $G$-invariant domain $\Omega$ yields reducing subspaces $\mathbb P_\varrho \mathcal H$ and $\mathbb P_\varrho^{ii}\mathcal H$, where $\varrho \in \widehat G, \, 1 \leq i \leq \deg \varrho$ and $\widehat G$ is the equivalence classes of irreducible representations of $G$ \cite{BDGS}. It follows from \cite[Proposition 2.2]{S R} and \cite[Proposition 1, p.556]{T} that the mapping
$
\boldsymbol\theta=(\theta_1,\ldots,\theta_n):\Omega\to\boldsymbol\theta(\Omega),
$
is a proper holomorphic map. In particular, $\boldsymbol \theta(\Omega)$ is a domain, and the quotient topological space $\Omega/G$ admits the structure of a complex analytic space biholomorphic to $\boldsymbol \theta(\Omega)$ - justifying the nomenclature  \textit{Quotient domains}. In \cite{S R}, it was shown that every proper holomorphic mapping which is factored by automorphism is intimately related to a basic polynomial map $\boldsymbol \theta$. Thus, generalizing from the symmetric group to finite complex reflection groups offers a new perspective on the problem of minimal reducing subspaces under proper holomorphic mappings  \cite{S DPW, Z} and naturally links to our works in \cite{BS, BGMS, BDGS}. Recently, this has been explored extensively in the literature; see, for instance, \cite{G, ZH}.

\emph{\textbf{On realizations on quotient domain.}} In \cite{BDGS}, it was established that the submodules $\mathbb P_\varrho \mathcal H$ admit a realization as analytic Hilbert modules of rank $({\rm deg}\,\varrho)^2$ over $\mathbb C[\boldsymbol z]$, consisting of holomorphic functions on the quotient domain $\boldsymbol \theta(\Omega)$. This result was obtained through a generalization of the Chevalley–Shephard–Todd Theorem to the setting of algebra of holomorphic functions. It plays a crucial role to obtain a substantial improvement of the earlier result \cite[Corollary 3.12]{BGMS} for the symmetric group $\mathfrak S_n$, where it was shown that for any partition $\boldsymbol p \vdash n$, the module $\mathbb P_{\boldsymbol p}\mathcal H$ is locally free of rank $({\rm deg}\,\boldsymbol p)^2$ only on an open subset of $\boldsymbol s(\Omega) \setminus \boldsymbol s(\mathcal Z)$. Here $\boldsymbol s:\mathbb C^n \to \mathbb C^n$ denotes the symmetrization map, and $\mathcal Z$ is the set of critical points of $\boldsymbol s$. This provides a powerful framework for developing geometric and analytic unitary invariants \cite{CD, CS}, offering an avenue towards  complete classification of the reducing submodules $\{\mathbb P_\varrho\mathcal H\}_{\varrho\in\widehat G}$ up to unitary equivalence.

In recent years, attention has broadened from the symmetrized polydisc to other quotient domains such as the tetrablock, monomial polyhedra, and fat Hartogs triangles. This line of investigation typically unfolds along two directions: on the one hand, extending results obtained for special domains to broader families of quotient domains, and on the other, applying general structure theorems to recover more refined information in concrete settings. For example, the Hilbert space method of \cite{MSZ} for computing Bergman kernels was extended in \cite{T} to arbitrary proper holomorphic maps, and in \cite{G1} to quotient domains which includes monomial polyhedra and fat Hartogs triangles as special cases. Parallel activity have taken place on Hardy spaces. In \cite{MSZ}, a natural notion of Hardy space on $\mathbb G_n$ was introduced and its reproducing kernel was obtained explicitly. This construction was later generalized in \cite{S GS} to quotients of bounded symmetric domains, where Brown–Halmos type theorems for Toeplitz operators were also established -  some of which were studied in the context of the symmetrized polydisc in \cite{S BDS, S DS}. Toeplitz operators on the Bergman space  were  investigated in \cite{GN}. More recently, the structure theory for pure isometries and invariant subspaces, in the spirit of the Beurling–Lax–Halmos theorem, was extended from $\mathfrak S_n$ to the broader family of reflection groups $G(m,p,n)$ \cite{BGNS}. The study of invariant subspaces is closely tied to the theory of inner functions, raising questions analogous to those addressed for the symmetrized bidisc in \cite{ALY0} on quotient domains. Along these lines, \cite{Bhowmik-Kumar} examined the structure of rational inner functions on $\boldsymbol\theta(\mathbb D^n)$ for $G=G(m,p,n).$  Notable progress has also been made in multiple directions: $L^p$-regularity of  Bergman and Szego projections on quotient domains \cite{BCEM, GG, HW}, the  Hankel operators \cite{BBS}, homogeneous analytic Hilbert modules under non-transitive action of automorphism group \cite{BDHKM}, among others. Taken together, these studies highlight the emergence of a rich and interconnected theory of Hilbert function spaces on quotient domains - one that intricately weaves together operator theory, complex geometry, and representation theory - thereby opening numerous avenues for further exploration.

\section{Trace formulas in operator theory}

\noindent \textit{By Arup Chattopadhyay [2003arupchattopadhyay@gmail.com] and \\ Saikat Giri [saikatgiri90@gmail.com] from IIT, Guwahati; and \\
Chandan Pradhan [chandan.pradhan2108@gmail.com] from IISc.}

\bigskip

Let $\mathcal{H}$ be a complex separable Hilbert space. For $1 \le p < \infty$, the Schatten class of order $p$ on $\mathcal{H}$ is denoted by $\mathcal{S}^{p}(\mathcal{H})$ \cite{SkTo_book}.
	
Trace formulas are central objects in operator theory, motivated both by physics and mathematics. In 1953, Krein \cite{Kr53}, building on earlier work of Lifshitz \cite{Li52}, established a fundamental formula for self-adjoint operator $H_{0}$ perturbed by a self-adjoint trace-class operator $V \in \mathcal{S}^{1}(\mathcal{H})$:\vspace*{-.07in}
\begin{align}\label{babytf}
\textnormal{Tr}\big(f(H_0+V)-f(H_0)\big)=\int_{\mathbb{R}} f'(\lambda)\,\xi_{H_0,V}(\lambda)\,d\lambda,
\end{align}
valid whenever the Fourier transform of $f'$ belongs to $L^{1}(\mathbb{R})$.
	
The function $\xi_{H_0,V}$, known as the spectral shift function (SSF), describes how the spectrum changes when the operator is perturbed. Since its introduction, the SSF has become a cornerstone of perturbation theory, spectral flow, mathematical physics, and noncommutative geometry \cite{AzCaDoSu09, AzCaSu07,SkTo_book}. In \cite{ChSi18}, Chattopadhyay--Sinha analyzed the classical Carey--Helton--Howe--Pincus trace formula via Krein's SSF. These advances inspired many generalizations of trace formulas to broader operator frameworks, including noncompact perturbations and non-normal operators.
	
Analogues of \eqref{babytf} for unitary operators and contractions have been extensively studied \cite{AlPe16,ChSi22,DySk14,Kr62,Kr87,MaNe15,MaNePe19}. Unlike previous work, Chattopadhyay--Sinha \cite{ChSi22} established a trace formula for pairs of contractions in which the SSF is supported on the unit disk, rather than the unit circle as in earlier results \cite{MaNe15,MaNePe19}.
	
Extending Krein's trace formula beyond trace-class perturbations has been difficult. For $V\in\mathcal{S}^{2}(\mathcal{H})$, Koplienko \cite{Ko84} replaced the left hand side of \eqref{babytf} with a Taylor remainder. Higher-order formulas were attempted in \cite{ChSi13,DySk09,Sk10}, but only for $\mathcal{S}^2(\mathcal{H})$-perturbations (see also \cite{Ne88a,Pe05,PoSu12}). A full generalization was established by Potapov--Skripka--Sukochev \cite{PoSkSu13} for $V\in\mathcal{S}^{n}(\mathcal{H})$, proving
\begin{align}\label{tf}
\textnormal{Tr}\Big(f(H_0+V)-\sum_{k=0}^{n-1}\frac{1}{k!}\frac{d^{k}}{dt^{k}}\Big|_{t=0}f(H_{0}+tV)\Big)=\int_{\mathbb{R}} f^{(n)}(\lambda)\,\xi_{n,H_0,V}(\lambda)\,d\lambda.
\end{align}
The proof relies on the machinery of multiple operator integrals (MOIs) \cite{SkTo_book}, a powerful tool in modern perturbation theory. MOIs were introduced by Peller \cite{Pe06} and, independently, by Azamov--Carey--Dodds--Sukochev \cite{AzCaDoSu09}, as a far-reaching generalization of Birman--Solomyak's classical double operator integral \cite{BiSo66}. Since their inception, MOIs have become a central technique in the analysis of operator perturbations, allowing for significant extensions of trace formulas such as \eqref{tf}. For instance, they have enabled versions of \eqref{tf} for contractions \cite{PoSkSu14}, for unitary operators \cite{PoSkSu16,Sk17}. In a related development, Aleksandrov--Peller \cite{AlPe11} established \eqref{tf} for functions in Besov spaces, extending the class of functions considered in \cite{PoSkSu13}. It is worth noting that recently, Coine--Le Merdy--Sukochev \cite{CoLeSu21} introduced a MOI that encompasses the largest class of scalar functions.
	
Explicit representations of SSFs are rare. To address this, Voiculescu \cite{Vo87}, Sinha--Mohapatra \cite{MoSi94,MoSi96}, and Potapov–Sukochev–Zanin \cite{PoSuZa14} employed finite-dimensional approximations. Later, Chattopadhyay--Sinha, followed by Chattopadhyay--Das--Pradhan, established second-order trace formulas for self-adjoint \cite{ChSi12}, unitary \cite{ChDaPr22}, and contraction operators \cite{ChDaPr24}. Higher-order cases, however, remain inaccessible.
	
Significant progress has been made on higher-order results for broader operator classes. Chattopadhyay--Pradhan--Skripka \cite{ChPrSk24} proved trace formulas for contractions and maximal dissipative operators, removing earlier restrictions \cite{ChPr23,MaMo06,MaNePe19,Sk17}, and showed that for contractions the higher-order SSF can be supported on the unit disk. For Dirac and Schr\"odinger operators with noncompact perturbations, formulas of type \eqref{tf} were obtained for relatively Schatten-class perturbations \cite{ChSk18,NuSk22,SkTeJOT,PoSkSu15,Sk14} and for $\tau$-compact resolvents \cite{ChPrSk23,SkAnJOT}, assuming bounded perturbations. More recently, Chattopadhyay--van Nuland--Pradhan \cite{ChNuPr24} treated unbounded perturbations: for self-adjoint $H_{0}$ and symmetric, relatively $H_{0}$-bounded $V$, they derived explicit formulas for $\frac{d^{n}}{dt^{n}}\big|_{t=0}f(H_{0}+tV)$, proved convergence of the noncommutative Taylor expansion, and established higher-order SSFs under natural summability assumptions. The last result extends the work of Yafaev \cite{Ya05} and Aleksandrov–Peller \cite{AlPe25} to the higher-order case.
	
Another direction concerns enlarging the admissible function class for \eqref{tf}. Peller \cite{Pe16} (self-adjoint case) and Aleksandrov–Peller \cite{AlPe16} (unitary case) identified optimal classes for Krein trace formulas, but none are known for higher-order cases. Recently, Chattopadhyay--Coine--Giri--Pradhan \cite{ChCoGiPr24b,ChCoGiPr24a} unified differentiability results for self-adjoint, unitary, and contractive operators, establishing second-order trace formulas for the broader scalar function class so far. In \cite{ChCoGiPr24b}, they further obtained modified higher-order trace formulas with optimal scalar function classes, extending earlier results of Skripka \cite{Sk17} and Bhattacharyya--Chattopadhyay--Giri--Pradhan \cite{BhChGiPr25}.
	
The situation is more subtle in several variables. Aleksandrov--Peller--Potapov \cite{AlPePo19} proved that no analogue of \eqref{babytf} holds for noncommuting self-adjoint tuples, and asked whether for commuting tuples of self-adjoint operators $\mathbf{H}_{n}=(H_{1},\dots,H_{n})$, $\mathbf{H}_{n}(1)=(H_{1}+V_{1},\dots,H_{n}+V_{n})$ with $V_{i}\in\mathcal{S}^{1}(\mathcal{H})$, there exist measures $\mu_{j}$ such that  \vspace*{-.08in}
\begin{align}\label{tfsev}
\textnormal{Tr}\big(f(\mathbf{H}_{n}(1))-f(\mathbf{H}_{n})\big)=\sum_{j=1}^{n}\int_{\mathbb{R}^{n}}\frac{\partial f}{\partial \lambda_{j}}\,d\mu_{j}.\vspace*{-.07in}
\end{align}
This is resolved assuming $\mathbf{H}_{n}(t):=(H_{1}+tV_{1},\ldots,H_{n}+tV_{n})$ commutes for all $t\in[0,1]$: first for bounded operators by Skripka \cite{Sk15}, then for unbounded ones by Chattopadhyay--Giri--Pradhan--Usachev \cite{ChGiPrUs23}. Moreover, \cite{ChGiPrUs23} shows that if $V_{i}$ lie in the Lorentz ideal and $\textnormal{Tr}$ is replaced by a bounded singular trace $\tau$ \cite{LoSuZa_book}, then \eqref{tfsev} holds without this additional commutativity assumption. Higher-order cases were studied in \cite{ChGiPr24}, but the general multivariable problem remains open. Notably, Chattopadhyay--Sinha \cite{ChSi14} studied a different type of trace formula (a Stokes-like formula) in the multivariable setting.
	
Thus, the theory of trace formulas--ranging from Krein's original result to modern higher-order and multivariable extensions--continues to evolve, revealing deep connections between perturbation theory, noncommutative analysis, and spectral theory.
	
Research on trace formulas and spectral shift functions is longstanding and active. Since providing a complete list of contributions is beyond the scope of this survey note, the authors may unfortunately omit some important references; however, further sources can be found in the works we cite.

\section{Dilation of operators and algebraic varieties}

\noindent \textit{By B. Krishna Das [bata436@gmail.com] from IIT Bombay.}

\bigskip

Dilation theory forms a cornerstone of operator theory, providing a framework to understand commuting contractions on Hilbert spaces via their extensions or dilations to isometries or unitaries. The foundational results by Sz.-Nagy and Foias (\cite{NF}), and its extensions by Ando (\cite{Ando}), provide dilation theorems for single and pairs of commuting contractions, respectively. However, beyond two variables, the existence of such dilations generally fails. The limitations of these classical results motivate further study of dilation theory and add the need to understand operator tuples that fall outside these known theorems. The following summary synthesizes advances made in several papers, jointly with my collaborations, on isometric dilations of tuples of commuting contractions and its connection to algebraic varieites.

The classical von Neumann inequality (\cite{vN}) provides a bound on the norm of a polynomial evaluated at a contraction via its supremum over the unit disc. While Ando’s theorem extends this result to pairs of commuting contractions, the resulting inequality is often not sharp. In pursuit of sharper estimates, Agler and McCarthy (\cite{AgMc}) showed that for a pair of strict commuting matrices, the von Neumann inequality holds on a \emph{distinguished variety}--an algebraic subset of the bidisc defined as an algebraic curve intersecting only the distinguished part of the boundary of the bidisc. Motivated by this, we investigated whether there is a more precise dilation result for commute pairs of Hilbert space operators, leading to a sharper form of the von Neumann inequality. For a contraction $T$, we say $T$ is pure if $T^{* n}\to 0$ in the strong operator topology, and we define the defect operator and the defect space of $T$ as $D_T:=(I-TT^*)^{\frac{1}{2}}$ and $\mathcal D_{T}:=\overline{\text{Ran}}(I-TT^*)^{\frac{1}{2}}$, respectively. Our main result in this context is the following explicit Ando-type dilation theorem and sharp von Neumann inequality.
\begin{theorem}[\cite{DS}]\label{Ando}
Let $(T_1, T_2)$ be a pair of commuting contractions such that $T_1$ is pure and the defect spaces of $T_1$ and $T_2$ are finite-dimensional. Then:
\begin{itemize}
\item[(a)] There exists a $\mathcal{B}(\mathcal D_{T_1})$-valued rational inner function $\Psi$ such that the pair $(M_z, M_{\Psi})$ of multiplication operators on $\mathcal D_{T_1}$-valued Hardy space dilates $(T_1, T_2)$.

\item[(b)] If $T_2$ is also pure then there exists a distinguished variety $\mathcal{V} \subset \overline{\mathbb{D}}^2$ such that
\[
\|p(T_1,T_2)\|\le \sup_{(z_1,z_2)\in\mathcal V}|p(z_1,z_2)|, \text{ for all }\ p\in\mathbb C[z_1,z_2].
\]

\end{itemize}
\end{theorem}
The inner function $\Psi$ that appears in this result can be explicitly realized through transfer function models, and the associated distinguished variety $\mathcal{V}$ admits a determinantal representation in terms of $\Psi$. This result extends and significantly strengthens the work of Agler and McCarthy by moving from matrices to Hilbert space contractions and offers a new, concise proof of their result. Moreover, the dilation result in part (a) of the theorem is \emph{optimal} in the following sense: if a pair of commuting contractions $(T_1, T_2)$ dilates to $(M_z, M_{\Psi})$ on an $\mathcal{E}$-valued Hardy space, where $\mathcal{E}$ is finite-dimensional and $\Psi$ is a $\mathcal{B}(\mathcal{E})$-valued rational inner function, then $T_1$ must be pure and the defect spaces of $T_1$ and $T_2$ must be finite-dimensional (see \cite{DS1}). The above result also illustrate that identifying operator tuples that dilate to specific classes of isometries can yield finer structural information about those tuples. This principle is further demonstrated by our dilation result in the context of BCL pairs of isometries, as described below. The celebrated result of Berger, Coburn, and Lebow (BCL) (\cite{BCL}) characterizes all possible ways in which the Hardy shift $M_z$ can be factorized as a product of commuting isometries. Identifying classes of pairs of commuting contractions that dilate to a pair of BCL isometries leads to the following factorization theorem for a contraction.
\begin{theorem}[\cite{DSS}]
Let $T$ be a pure contraction and let $T\cong P_{\mathcal Q}M_z|_{\mathcal Q}$ be the Sz.-Nagy and Foias representation of $T$. Then:
\begin{itemize}
\item[(a)] A pair of commuting contractions $(T_1, T_2)$ with $T_1T_2$ is pure dilates to a BCL pair of isometries
on some vector-valued Hardy space.
\item[(b)] If $T=T_1T_2$ for some pair of commuting contractions $(T_1, T_2)$, then there exist $\mathcal{B}(\mathcal D_T)$-valued polynomials $\varphi$ and $\psi$ of degree one such that
\[
(T_1, T_2)\cong (P_{\mathcal Q}M_{\varphi}|_{\mathcal Q}, P_{\mathcal Q}M_{\psi}|_{\mathcal Q}).
\]
\end{itemize}
\end{theorem}
The dilation result in part (a) has been applied to derive several other dilation theorems. For example, it yields an alternative proof of Ando’s theorem (\cite{MB}) and has been used to construct dilations for various classes of $n$-tuples ($n \ge 3$) of commuting contractions defined by certain positivity conditions (see \cite{BD, BDHS, BDS, BDDP}).

A polynomial in two variables whose zero set in the closed bidisc is a distinguished variety is called a \emph{toral polynomial}. A pair of commuting contractions is called a \emph{toral pair} if it is annihilated by a toral polynomial. A consequence of Theorem~\ref{Ando} is that if $(T_1, T_2)$ is a pair of commuting contractions such that both $T_1$ and $T_2$ are pure contractions with finite defect, then $(T_1, T_2)$ is toral and admits an isometric lift that is also toral. This naturally leads to the following constrained Ando dilation problem, which is closely related to the rational dilation problem on a distinguished variety of the bidisc: \emph{Does a toral pair of commuting contractions admit a toral isometric lift?} In general, the answer is negative. However, we have the following positive result.
\begin{theorem}[\cite{DS2}]
Any toral pair of commuting contractions with finite defect spaces admits a toral isometric lift.
\end{theorem}
This theorem implies that if a pair of commuting contractions with finite defect spaces has a distinguished variety of the bidisc as its spectral set, then there exists (possibly a different) distinguished variety that serves as a complete spectral set. It is now well known that the rational dilation problem on a distinguished variety does not hold in general; the Neil parabola is one such example where rational dilation fails. One of the questions we aim to address is the following:

\textbf{Question:} \emph{Under what conditions does rational dilation hold on a distinguished variety of the bidisc?}

In conclusion, these works significantly advance the understanding of operators through dilation theory and sharp von Neumann inequalities on distinguished varieties. A key feature in each of these works has been the construction of explicit dilation. Looking ahead, these advances also point toward several challenging open problems, most notably the rational dilation problem on distinguished varieties of the bidisc.

\section{The Weighted Discrete Semigroup Algebra $\ell^1(S, \omega)$}

\noindent \textit{By H. V. Dedania [hvdedania@gmail.com] from Sardar Patel Univeristy.}

\bigskip

A \emph{Banach algebra} is an associative algebra $A$ together with a norm $\|\cdot\|$ such that $(A, \|\cdot\|)$ is a Banach space and $\|ab\| \leq \|a\|\|b\| \; (a, b \in A)$. The classical examples of Banach algebras are $C_0(\Omega)$, $BL(X)$, and $L^1(G)$; where $\Omega$ is a locally compact, Hausdorff space, $X$ is a Banach space, and $G$ is a locally compact, Hausdorff topological group. Many important algebras in analysis and in other branches of mathematics are Banach algebras.  There is a strong bonding between the algebraic structure and the topological structure in Banach algebras. The algebraic structure of $A$ determines various properties of its topological structure. For example, each multiplicative linear functional on $A$ is continuous \cite[Theorem 2.1.29]{Dales}.

Let $A$ be a Banach algebra and $A^*$ be its Banach space dual. Let $S(A)= \{x \in A: \|x\|=1\}$, $D_{A}(x) = \{\varphi \in A^*: \|\varphi\|=1=\varphi(x)\}$, and $V_{A}(a;x) = \{\varphi(ax): \varphi \in D_{A}(x)\}$ for $a \in A$ and $x \in S(A)$. Then the set $V_{A}(a) =   \cup\{V_{A}(a;x): x \in S(A)\}$ is called the \emph{spatial numerical range (SNR)} of $a$ in $A$ with respect to $\|\cdot\|$. The SNR $V_{A}(a)$  highly depends on both the algebra $A$ and the norm $\|\cdot\|$. It is proved in \cite[Theorem 2.1]{Dedania-2} that $V_{A}(a;x)$ is always convex. So the following is a natural question: Is the SNR $V_{A}(a)$ always convex? We do not know its answer. If $A$ has an identity $e$ with $\|e\|=1$, then $V_{A}(a) = V_{A}(a; e)$ is convex \cite[Corollary 2.1]{Dedania-2}. It is convex in many non-unital Banach algebras \cite{Dedania-3,Patel}. On the other hand, Bonsall and Duncan defined the SNR for any bounded linear operator $T \in BL(X)$, which is slightly different from our definition. They constructed an operator $T$ on ${\mathbb C}^2$ whose SNR $V_{BL(X)}(T)$ is not convex \cite[Page-98]{BonsallAndDuncan}. Unfortunately, this operator $T$ does not help us to find an element $a$ in $A$ such that $V_A(a)$ is not convex.

One of the most important Banach algebra is the weighted discrete semigroup algebra $\ell^1(S, \omega)$, where $S$ is a semigroup and $\omega$ is a \emph{weight} on $S$, i.e., $0 < \omega(s)$ and $\omega(st) \leq \omega(s)\omega(t)$ for all $s, t \in S$ \cite[Page-159]{Dales}. The Banach algebra norm on it is denoted by $\|\cdot\|_{\omega}$. The Banach algebra $\ell^1(S, \omega)$ is mostly used as a counter example. However, it is studied independently also. Its Banach algebra structure is highly influenced by $S$ and $\omega$. For example, it is commutative iff $S$ is abelian; and it is semisimple iff $S$ is separating and $\omega$ is semisimple \cite[Theorem 2.3]{Dedania-5}. Here we present some results on $\ell^1(S, \omega)$, which are proved by us.

If $A$ is a semisimple Banach algebra, then we get another norm on $A$; namely, the \emph{operator norm} $\|\cdot\|_{op}$, which is defined as $\|a\|_{op} = \sup\{\|a x\|: x \in A, \|x\| \leq 1\} \; (a \in A)$. In general, $\|\cdot\|_{op} \leq \|\cdot\|$ on $A$. The norm $\|\cdot\|$ is \emph{regular} if $\|\cdot\|_{op} = \|\cdot\|$ on $A$. For example, every $C^{\ast}$-norm is a regular norm. A natural question is the following: When is the weighted norm $\|\cdot\|_{\omega}$ regular?, i.e., $\|\cdot\|_{\omega op} = \|\cdot\|_{\omega}$ on $\ell^1(S, \omega)$? We have proved the following.\\
\textbf{Theorem-(A):} \cite{Dedania-4} Let $S$ be a right cancellative semigroup and $\omega$ be a weight on $S$. Let $k \in {\mathbb N}$ and $\widetilde{\omega}_{0} = \omega$. Define $\widetilde{\omega}_k(s)= \sup \{\frac{\widetilde{\omega}_{k-1}(st)}{\widetilde{\omega}_{k-1}(t)} : t \in S\} \; (s \in S)$. Then
\begin{enumerate}
  \item Each $\widetilde{\omega}_k$ is a weight on $S$.
  \item $\widetilde{\omega}_{k+1} \leq \widetilde{\omega}_{k}$ on $S$ for each $k$.
  \item $\ell^1(S, \omega) \subseteq \ell^1(S, \widetilde{\omega}_{k-1}) \subseteq \ell^1(S, \widetilde{\omega}_{k})$.
  \item $(\ell^1(S, \widetilde{\omega}_{k}), \, \|\cdot\|_{\widetilde{\omega}_{k}})$ is a Banach algebra.
  \item $\lim\limits_{n \rightarrow \infty} \widetilde{\omega}_k(s^n)^{\frac{1}{n}} = \lim\limits_{n \rightarrow \infty} \widetilde{\omega}_{k+1}(s^n)^{\frac{1}{n}} \; (s \in S)$.
  \item $\|\delta_t\|_{\widetilde{\omega}_{k} op} = \|\delta_t\|_{\widetilde{\omega}_{k+1}} \; (t \in S)$.
  \item $\|f\|_{\widetilde{\omega}_k op} \leq \|f\|_{\widetilde{\omega}_{k+1}} \; (f \in \ell^1(S,\omega))$.
  \item $\widetilde{\omega}_{k-1}$ has \emph{F-property} (i.e., for each finite set $\{t_1, \ldots, t_n\} \subset S$ and $r < 1$, there is $s \in S$ such that $\frac{\widetilde{\omega}_{k-1}(t_i s)}{\widetilde{\omega}_{k-1}(s)} \geq r \widetilde{\omega}_k(t_i) \; \; (1 \leq i \leq n)$) iff $\|\cdot\|_{\widetilde{\omega}_k op} = \|\cdot\|_{\widetilde{\omega}_{k+1}}$ on $\ell^1(S, \omega)$.
  \item The $1$-norm $\|\cdot\|_{1}$ is regular on $\ell^1(S)$.
\end{enumerate}
On the other hand, for $1 \leq p < \infty$, the $p$-norm $\|\cdot\|_{p}$ is never regular on the sequence algebra $\ell^p$ with pointwise multiplication \cite[Theorem 2.8]{Dedania-4}.

An element $a$ in a commutative Banach algebra $A$ is said to be \emph{compact} if the linear operator $L_a : A \longrightarrow A; \, x \longmapsto ax$  is a compact operator in the norm topology. Let $K(A)$ denote the set of all  compact elements in $A$. Then $K(A)$ is a closed ideal in $A$. If $A$ has finite dimension, then clearly $K(A) = A$. On the other hand, if the Gel'fand space $\bigtriangleup(A)$ of a semisimple, commutative Banach algebra $A$ has no isolated point, then $K(A) = \{0\}$. Thus, the concept of compact elements is interesting mainly in infinite dimensional, radical, commutative Banach algebras; the convolution algebras $L^1({\mathbb R}_+, \omega)$ and $l^1(S, \omega)$ are such Banach algebras.

The set of all compact elements of $L^1({\mathbb R}_+, \omega)$ has been studied by Bade and Dales and for $l^1(S, \omega)$ by Gr{\o}nb{\ae}k. Some more results on compact elements in $l^1(S, \omega)$ have been studied by the author. We recently proved some weighted discrete analogues of the results proved for $L^1({\mathbb R}_+, \omega)$ about compact elements \cite{Dedania-1}. We also set the record right by correcting some statements claimed by Gr{\o}nb{\ae}k. Let $K(S,\omega)$ be the set of all compact elements in $\ell^1(S,\omega)$. Our main results are the following two.\\
\textbf{Theorem-(B):} \cite{Dedania-1} Let $\omega$ be a weight on an abelian semigroup $S$ and $f \in \ell^1(S,\omega)$. Consider the following three statements:\\
(a) $\sum\{|f(s)| \omega_{0}(s) : s \in S\} = 0$; \; (b) $\lim\limits_{s\rightarrow [\infty]} \widetilde{\delta}_s \ast f =0$; \; (c) $f \in K(S,\omega)$.\\
Then
\begin{enumerate}
\item $(a) \Longrightarrow (b) \Longrightarrow (c)$.
\item If $S$ is weakly cancellative, then $(c) \Longrightarrow (b)$.
\item If $S$ is cancellative, then $(b) \Longrightarrow (a)$.
\end{enumerate}
\textbf{Theorem-(C):} \cite{Dedania-1} Let $S$ be cancellative and abelian, and $\omega$ be a weight on $S$. Then there is a semigroup ideal $T$ in $S$ such that $K(S, \omega) = \{f \in \ell^1(S, \omega) : \textrm{supp}f \subseteq T\}$.\\

\noindent
\textbf{Open Problems:} Following are some open problems in this area.
\begin{enumerate}
\item Let $S$ be an infinite semigroup. Is there always an unbounded weight on $S$?
\item Let $G$ be a discrete group. Is $\ell^1(G,\omega)$ always semisimple?
\item Is there a Banach algebra $A$ such that $V_A(a)$ is not convex for some $a \in A$?
\item Is there an algebra $A$, which has finitely many (but more than one) non-equivalent Banach algebra norms?
\item When does the converse of Theorem-(C) hold true?, i.e., If $T$ is a semigroup ideal in $S$, then does there exist a weight $\omega$ on $S$ such that $K(S, \omega) = \ell_T^1(S, \omega)$?
\end{enumerate}

\section{Cartan isometries and Cartan contractions}

\noindent \textit{By Surjit Kumar [surjit@iitm.ac.in] from IIT Madras; and Paramita Pramanick [paramitapramanick@gmail.com] from ISI Kolkata.}

\bigskip

Let $\mathcal H $ be a complex separable  Hilbert space. An operator $T$ on $\mathcal H$
 is an isometry if $T^*T=I.$ A canonical example of isometry is the unilateral shift on $\ell^2(\mathbb{N})$.
By a commuting $d$-tuple $\boldsymbol T,$ we mean a tuple of operators $T_1,\cdots, T_d $ in $ \mathcal B(\mathcal H)$ that commutes pairwise. A {\it spherical isometry} is a commuting $d$-tuple $\boldsymbol{T}$ which satisfies the condition $ T_1^* T_1+ \cdots +T_d^* T_d=I.$ A remarkable result of Athavale shows that every spherical isometry is subnormal
with normal spectrum contained in the unit sphere $\partial \mathbb B$ (see \cite{At1}). Consequently, one obtains a natural functional calculus for the ball algebra. Thus, a spherical isometry may be regarded as a multi-variable analogue of an isometry associated with the unit ball.
This note explores isometries arising from classical Cartan domains, referred as Cartan isometries and associated contractions.

Let $\Omega$ be a classical Cartan domain of rank $r$ in $\mathbb{C}^d$ of type $(r, a, b).$ Let $S_{\Omega}$ be the Shilov boundary of the domain $\Omega.$ Let $G$ be the connected component of the identity in $\rm{Aut}(\Omega)$, the biholomorphic automorphism group of $\Omega$. Let $\mathbb{K}=\{g \in G : g(0)=0\}$ be a maximal compact subgroup of $G$.  Note that $\mathbb K$ is a subgroup of linear automorphism in $G.$
 The group $\mathbb K$ acts on the space of analytic polynomials $\mathcal P$ by composition. Under this action, $\mathcal P$ decomposes into irreducible, mutually $\mathbb{K}$-inequivalent subspaces $\mathcal P_{\underline s}$ such that $\mathcal{P}=\sum \mathcal{P}_{\underline s },$ where the sum runs over all signatures, $\underline s=(s_1, \ldots , s_r)\in \mathbb Z^r_+,\; s_1 \geq \ldots \geq s_r\geq 0$. Let $\{\psi_\alpha^{\underline s}(z)\}_{\alpha=1}^{d_{\underline s}}$ be an orthonormal basis of $\mathcal P_{\underline s}$ with respect to the Fischer-Fock inner product, where $d_{\underline s}$ is the dimension of $\mathcal P_{\underline s}.$ The reproducing kernel $K_{\underline s}$ is given by
$K_{\underline s}(z,w)= \sum_{\alpha=1}^{d_{\underline s}} \psi_\alpha^{\underline s}(z) \overline{\psi_\alpha^{\underline s}(w)}.$

A Cartan isometry is a commuting subnormal $d$-tuple $\boldsymbol{T}$ with normal spectrum $\sigma_n(\boldsymbol{T})$ contained in $S_{\Omega}$ (see \cite{At}).
A well-known example of a Cartan isometry is the Szeg\"o shift on the Hardy space $H^2(S_{\Omega}).$ In fact, a $d$-tuple $\boldsymbol T \in \mathcal AK(\Omega)$ is a Cartan isometry if and only if $\boldsymbol T$ is unitarily equivalent to the Szeg\"o shift. For more details on $\mathcal AK(\Omega),$ a certain class of $\mathbb K$-homogeneous operator tuples, we refer to \cite{GKP}.
In \cite{At}, Athavale provided a characterization of Cartan isometries addressing each classical Cartan domain individually.
Recently, in \cite{KMP2025}, we have characterized Cartan isometries for all types of classical Cartan domains, unifying them under a single framework.
\begin{theorem}\label{characterization}
A commuting $d$-tuple $\boldsymbol{T}$ is a Cartan isometry if and only if  \begin{equation} \label{Cartaniso}\Delta^{(\ell)}(z,w)(\boldsymbol{T},\boldsymbol{T}^*)=\binom{r}{\ell} I_{\mathcal{H}}\; \mbox{~for~all~} 1 \leq \ell \leq r,\end{equation}
where $(\ell):= (1,\ldots,1,0,\ldots,0)$ denotes  the signature with first $\ell$ many ones and
\begin{equation} \Delta^{(\ell)} (z,w) = \prod_{j=1}^{\ell} (1+\frac{a}{2}(j-1))K_{(\ell)}(z,w).\end{equation}
\end{theorem}
Consequently, we have the following results:
\begin{corollary}
 Let $\boldsymbol T$ be a Cartan isometry. The following holds true:
 \begin{enumerate}
     \item [(i)] The Taylor spectrum $\sigma(\boldsymbol T)$ is contained in $\overline{\Omega}.$
     \item [(ii)] For every $g \in G, \; g(\boldsymbol T)$ is a Cartan isometry.
     \item [(iii)] Every Cartan isometry is reflexive.
 \end{enumerate}
\end{corollary}
The Wold–von Neumann decomposition \cite[ pp. 14]{Co} asserts that every isometry on a Hilbert space is unitarily equivalent to the direct sum of a unilateral shift (of possibly countable multiplicity) and a unitary operator. This result plays a central role in Beurling’s characterization of $z$-invariant subspaces of the Hardy space over the unit disc \cite[pp.23 ]{Co}.
Using the model theory of {\it row $d$-hypercontractions}, Athavale showed that a spherical isometry $\boldsymbol{T}$ admits a Wold-type decomposition exactly when $\boldsymbol{T}^*$ is a $d$-hypercontraction \cite[Theorem 2.1]{At3}. This naturally leads to the study of a Wold-type decomposition for Cartan isometries.

Recall that a commuting $d$-tuple $\boldsymbol T$ is a {\it spherical contraction} if $ T_1^* T_1+ \cdots +T_d^* T_d \leq I.$ Equivalently, the operator $D_{\boldsymbol T}: \mathcal H \to \mathcal H^{(d)}$ given by $D_{\boldsymbol T}h=(T_1h, \ldots, T_dh)$ is a contraction. Analogously, we define a notion of contraction for classical Cartan domains.
For $m\geq n\geq 1$, let $\Omega_{n,m}=\{z\in \mathcal M_{n\times m}: \|z\|<1\}$ be the Cartan domain of type-I.
Let $\boldsymbol T=(T_{11},\ldots,T_{1m},T_{21},\ldots,T_{2m},\ldots,T_{n1},\ldots,T_{nm})$ be a commuting $mn$-tuple of operators on $\mathcal H$. Consider $D_{\boldsymbol T}:\mathcal H^{(n)}\to \mathcal H^{(m)}$ given by
$D_{\boldsymbol T}(h_1,\ldots,h_n)=(\sum_{i=1}^nT_{i1}h_i,\ldots,\sum_{i=1}^nT_{im}h_i)$ for all $h_i\in \mathcal H$ and $i=1,\ldots,n.$ The $mn$-tuple $\boldsymbol T$ is said to be a {\it Cartan contraction of type}-I if  $D_{\boldsymbol T}$ is a contraction. A routine verification shows that $\boldsymbol T$ is a Cartan isometry if and only if  $D_{\boldsymbol T}$ is an isometry. In a similar fashion, one can define Cartan contractions for Cartan domains of type-II and type-III. A realization formula for Schur-Agler class associated with such (strict) contraction appeared in the work of \cite{BB2004}.  If $\boldsymbol T$ is a Cartan contraction of type-I, then by \cite[Lemma 1]{AT2003}, $\sigma(\boldsymbol T)\subseteq \overline{ \Omega}_{n,m}.$ Also,  $g(\boldsymbol T)$ is a Cartan contraction of type-I for all $g \in G$.

Consider the $mn$-tuple $\boldsymbol M^{(\nu)}$ of multiplication operators given by $M_{ij}(f)(z)=z_{ij}f(z)$ on the weighted Bergman space $ \mathcal H^{(\nu)}$ of holomorphic function defined on $\Omega_{n,m}.$ Then, for any $\nu \in \{ m, \ldots, m+n-1\}\cup(m+n-1, \infty), \; \boldsymbol M^{(\nu)}$ is a subnormal $d$-tuple, and therefore, is a Cartan contraction of type-I. Furthermore, for $\nu< m, \; \boldsymbol{M}^{(\nu)}$ is not a Cartan contraction of type-I.

A fundamental result due to von Neumann \cite{von} states that if $T$ is a contraction on a Hilbert space $\mathcal H$, then $\|p(T)\| \leq \|p\|_{\infty, \mathbb D}:=\sup \{|p(z)|: |z|<1\}$ for every polynomial $p$. In other words, the induced homomorphism $\rho_T(p)=p(T)$ is contractive whenever $T$ is a contraction. Sz.-Nagy's dilation theorem \cite{Sz} shows that $\rho_T$ is not only contractive but {\it completely contractive}, that is, $T$ admits a unitary (power) dilation. Equivalently, $\|P(T)\| \leq \|P\|_{\infty, \mathbb D}$ for every $n \times n$ matrix valued polynomials $P(z)=(p_{i,j}(z))$. Ando's dilation theorem \cite{Ando} extends the result to pairs of commuting contractions. Parrott \cite{Parrott1970} constructed an example of three commuting contractions that satisfies von Neumann’s inequality but fails to admit a commuting unitary dilation.
We say that a commuting $d$-tuple $\boldsymbol T$ admits a {\it spherical dilation} (normal $\partial \mathbb B$-dilation in Arveson’s terminology) if there exist a Hilbert space $\mathcal K \supset \mathcal H$ and a spherical isometry $\boldsymbol S=(S_1,\ldots, S_d)$ in $\mathcal B(\mathcal K)$ such that $\boldsymbol T^\alpha= P_{\mathcal H}\boldsymbol S^\alpha|{\mathcal H}$ for all $\alpha \in \mathbb Z^d_+$, where $P_{\mathcal H}$ denotes the orthogonal projection of $\mathcal K$ onto $\mathcal H$.  Due to a deep result of Arveson \cite{S Ar}, it follows that $\boldsymbol T$ admits a normal $\partial \mathbb B$-dilation if and only if the induced homomorphism $\rho_{\boldsymbol T}$ is completely contractive.

In \cite{At0}, Athavale established that every $d$-hypercontraction admits a spherical dilation (see also \cite{M-V}). Identifying an appropriate analogue of a $d$-hypercontraction in the setting of Cartan domains, which could potentially yield a Cartan isometric dilation, appears to be a challenging problem. On the other hand, for any spherical contractive  $d$-tuple  $\boldsymbol M$ of multiplication operators by the coordinate functions on a reproducing kernel Hilbert space with a $\mathcal U(d)$-invariant kernel, the induced homomorphism $\rho_{\boldsymbol M}$ is contractive, that is, $\|P(\boldsymbol M)\| \leq \|P\|_{\infty, \mathbb B}$ (cf. \cite[Proposition 2.5]{Ku}). In fact, by employing slice representation techniques (see \cite[Theorem 1.10]{Ch-K}), one can obtain a spherical dilation for such operator tuples. Looking ahead, it is natural to expect that the techniques of ``slice representation" for $\mathbb K$-invariant kernels may eventually lead to a Cartan isometric dilation for Cartan contractive multiplication tuples. Equivalently, this would establish that for any Cartan contractive $d$-tuple in $\mathcal{A}\mathcal{K}(\Omega),$ the induced homomorphism is completely contractive.

\section{Norms on operator space tensor products and $C^*$-envelopes of operator systems}

\noindent \textit{By Ajay Kumar [ak7028581@gmail.com] from University of Delhi; and Preeti Luthra [preeti@ms.du.ac.in] from Mata Sundri College for Women, University of Delhi.}

\bigskip

Von Neumann in collaboration with Murray began a mathematical quantization program. The key idea being to replace function spaces by $*$-algebra of bounded operators on Hilbert spaces. An {\em operator space} is a Banach space $X$ with an extra matricial norm structure. We have norms on $M_n(X)$ of $n \times n$ matrices with entries from $X$, where these norms must satisfy certain consistency requirements. The systematic study of operator spaces has been done mainly by Effros and Ruan \cite{effros, ruan}, Blecher and Paulsen \cite{Blecherp}, Blecher and Smith \cite{BlecherS} and Haagerup and Musat \cite{HaagM}.

Grothendeick, in his celebrated paper \cite{Groth}, considered $14$ natural norms defined on tensor product of two normed spaces. The foundation for a systematic study of tensor products of subspaces of $C^*$-algebras was done in \cite{effros,Blecherp}. They defined various operator space tensor norms, namely, the Haagerup norm, the operator space projective and injective norms.

The Haagerup norm $\|\cdot\|_h$ and Banach space projective norm $\|\cdot\|_{\gamma}$ are equivalent on the algebraic tensor product $A \otimes B$ of two $C^*$-algebras $A$ and $B$ if and only if either $A$ or $B$ is finite dimensional or $A$ and $B$ are infinite dimensional and subhomogeneous. The same conclusion holds if the Banach space projective norm is replaced by the operator space projective norm. The equivalence of the Banach space projective norm and operator space projective norm $\|\cdot\|_{\wedge}$ holds if and only if $A $ or $B$ are subhomogeneous \cite{kumarsinc}. The Haagerup and Banach space projective norm may be equivalent for an infinite-dimensional ternary ring of operator and an infinite dimensional $C^*$-algebra which are not subhomogeneous. No clear criterion is available for arbitrary operator spaces. Moreover, equivalence of Schur norm $\|\cdot\|_s$ and $\|\cdot\|_{\wedge}$ is still not known even for $C^*$-algebras \cite{itoh,kumar}. Basic tools involved are Grothendieck inequality \cite{haag,HaagM,PisierS}, lifting maps to second duals and to tensor product of second duals.

If $A$ and $B$ are $C^*$-algebras and $\|\cdot\|_{\alpha}$ and $\|\cdot\|_{\beta}$ are tensor norms on $A \otimes B$ such that $\|\cdot\|_{\alpha} \leq \|\cdot\|_{\beta}$, the identity map $i: A \otimes_\beta B \rightarrow A \otimes_\alpha B$ which maps elementary tensor into itself is continuous, and so can be extended uniquely to $i: A \otimes^\beta B \rightarrow A \otimes^\alpha B $. Injectivity when $\beta =\gamma$ and $\alpha=\lambda$ (the Banach space injective tensor norm) and $\beta=\wedge$ and $\alpha=\lambda$ were proved in \cite{haag,JainK}. We still do not know the injectivity in case $\beta =s$ or when $A$ and $B$ are operator spaces \cite{itoh}.
The Haagerup tensor product $M \otimes^h B$, where $M$ is a ternary ring of operators (TRO) and $B$ a $C^*$-algebra has been investigated. In contrast to the $C^*$-algebra, a counterexample can be obtained showing that the Haagerup and maximal tensor norms can be equivalent even when both $M$ and $B$ are infinite dimensional using inductive limits \cite{AntonyKL,KansalK}. We also established an isometric embedding of $M^{**} \otimes^h B^{**}  $ into $(M \otimes^h B)^{**}$. Although $M\otimes^h B$ does not admit a natural $C^*$-algebra or ternary ring structure, in fact we proved that $A \otimes^h B$ is a $C^*$-tring if and only if $A =\mathbb{C}$ or $B=\mathbb{C}$, thereby clarifying the boundary of when such structure survives even for $C^*$-algebras. Nevertheless, we developed an ideal-theoretic framework, classifying maximal, prime, and primitive ideals directly in terms of the ideal of $M$ and $B$. These results illustrate that the ideal structure of $M \otimes^h B$ retains much of the algebraic clarity seen in the $C^*$-algebraic context \cite{KansalK,KansalK2}.

Beside $C^*$-algebras and operator spaces, operator systems have attracted considerable attention of operator algebraists in recent years. \emph{Operator systems} are self-adjoint unital subspaces of $B(H)$ for a complex Hilbert space $H$ \cite{Paulsen}. Naturally, they are matrix ordered and matrix normed spaces just as $C^*$-algebras. Arveson in \cite{Arveson} remarked that one can generate non-$*$-isomorphic $C^*$-algebras from an operator system $S$ depending upon the representation, these are called as $C^*$-covers of $S$. The existence of minimal amongst all these $C^*$-covers for an operator system $S$, referred to as the $C^*$-envelope of $S$ and denoted by $C^*_e(S)$, was shown by Hamana \cite{hamana}.

In 2011, a new approach to the study operator system tensor products was introduced and a lattice structure of $6$ operator system tensor products was provided in \cite{KPTT1}. $\mathrm{ess}$ tensor product related to $C^*$-envelopes was later introduced in \cite{disgrp}, via the embedding $S \otimes_{\mathrm{ess}} {T} \subset C^*_e(S) \otimes_{\mathrm{max}} C^*_e({T})$. The whole theory of tensor products, even the very recent one, revolves around the work of Lance \cite{lan}. It is well known that the operator space minimal (injective) tensor product and maximal (projective) tensor product do not coincide even for finite dimensional matrix algebras, e.g., $\|E_{1j} \otimes E_{jj}\|_{\vee}= 1 \neq \sqrt{n} =\|E_{1j} \otimes E_{jj}\|_{\wedge}$ (see \cite{kumarsinc}), so the term nuclearity can not be extended to operator spaces. But due to the introduction of several tensor products in the category of operator systems, the notion of nuclearity was generalized to operator systems from $C^*$-algebras in \cite{KPTT1}.

The relationship of nuclearity of an operator system with the nuclearity of its $C^*$-envelope was considered in \cite{GL}. It was shown that nuclearity of the $C^*$-envelope of an operator system is equivalent to $(\min, \mathrm{ess})$-nuclearity of operator system. Moreover, a unital $C^*$-algebra is $(\min,\mathrm{ess})$-nuclear as an operator system if and only if it is nuclear as a $C^*$-algebra \cite{GL}. Associated to the the group $C^*$-algebra are finitely generated operator systems called \emph{group operator systems} \cite{disgrp}. The exhaustive list of nuclear group operator systems associated to minimal generating sets of finitely generated groups was provided in \cite{GL}. Unlike $C^*$-algebras, separable exact operator systems need not embed into the Cuntz algebra $\mathcal{O}_2$ \cite{Kirch,lupini}. It is the exactness of the $C^*$-envelope, rather than that of the operator system that makes an operator system embeddable into $\mathcal{O}_2$ \cite{LK1}. For a separable operator system $S$, the $C^*$-envelope $C^*_e(S)$ is exact if and only if there exist a unital complete order embedding of $S $ into $\mathcal{O}_2$. The commutativity of inductive limits and $C^*$-envelopes was studied in \cite{LK2}.

The study on the quantitative theory of polynomials in operator spaces was initiated in \cite{defant}. The $\lambda$-theory, based on the tensor norms obtained from homogeneous polynomial, was extended to matrix ordered spaces and Banach $*$-algebras in \cite{LKR}.

Dosiev in \cite{dosiev} proposed a concrete structure of local operator systems using the multinormed $C^*$-algebra $C^*_\mathcal{E} (\mathcal{D})$. Tensor product structure in the category of  local operator systems was introduced in \cite{beniwalKL1} and nuclearity in this category was discussed in \cite{beniwalKL3}. Further Arveson's notion of hyperrigidity of operator systems was extended to local operator systems in \cite{beniwalKL2}.

As a final remark, it is important to highlight that some recent studies e.g.~\cite{Cleve, Kim, Kennedy} have integrated the theoretical framework of operator algebras with practical explorations in quantum information theory, thereby emphasizing the significance of their deep and far-reaching interconnections.

\section{Bohr's Inequality and the Gleason--Kahane--\.Zelazko Theorem}

\noindent \textit{By Sneh Lata [sneh.lata@snu.edu.in] from Shiv Nadar Institution of Eminence; and Dinesh Singh [dineshsingh1@gmail.com] from O. P. Jindal Global university.}
\bigskip

Sneh Lata and Dinesh Singh, along with their collaborators, have made contributions–deemed significant by many–to the fields of complex analysis and functional analysis, particularly in two notable areas: Bohr's inequality and the Gleason--Kahane--\.Zelazko (GKZ) theorem. Their work extends important classical results that have strong implications for modern analysis to more general settings, including non-commutative spaces and new classes of function spaces.

\vspace{.3 cm}

\noindent {\bf Bohr's Inequality}

\vspace{.2 cm}

The classical Bohr's inequality, established by Harald Bohr \cite{Bohr} in number theory, states that for $f(z)=\sum_{k=0}^{\infty}a_{k}z^{k}\in H^{\infty}(\mathbb{D}),$
the inequality
$$
\sum_{k=0}^{\infty}|a_{k}|r^{k}\le \|f\|_{\infty}, \quad 0\le r\le \tfrac{1}{3},
$$
holds. The bound $r_0 = 1/3$ is optimal, meaning the inequality fails for some functions when $r > 1/3$. Bohr proved it initially for $r\le 1/6$, and later Riesz \cite{Dixon},
Schur \cite{Dixon}, Sidon \cite{Sidon}, Wiener \cite{Bohr}, and Tomic \cite{Tomic} extended it to $r=1/3$. More recently, Lata and Singh established analogues in non-commutative Hardy spaces associated with semifinite von Neumann algebras, yielding versions for von Neumann-Schatten class operators $\mathcal{C}_1$ and finite matrices. Earlier, Singh with Paulsen \cite{PS} obtained a general form for Hardy spaces associated with uniform algebras, showing the inequality’s independence from Fourier coefficients and its validity for more general coefficients. See also \cite{PPS, PS2}. Lata and Singh \cite{LS} later extended this analogue to non-commutative Hardy spaces associated with a von Neumann algebra $\mathcal{M}$ equipped with a normal semifinite faithful trace $\tau$.

\begin{theorem}\cite[Theorem 4.7]{LS}(Bohr's Inequality for Non-Commutative Hardy Spaces)\label{bohr1}
Let $x \in H^1(\mathcal{M})$ such that $\tau(x)\ge 0$ and $Re(x)\le y$ for some self-adjoint $y\in H^1(\mathcal{M})$ with finite trace. For a sequence
$\{x_m\}_{m=1}^{\infty}$ in $H_0^{\infty}(\mathcal{M})$ with $||x_m||_{\infty}\le1$, let $\alpha_0=\tau(x)$ and $\alpha_m=\tau(x x_m^{*})$ for $m\ge1$. Then:
\begin{equation}\label{bohr2}
\sum_{m=0}^{\infty}|\alpha_{m}|r^{m}\le\tau(y)
\end{equation}
whenever $0\le r\le\frac{1}{3}.$
\end{theorem}

This generalized theorem yields the classical Bohr's inequality when $\mathcal{M}=L^{\infty}(\mathbb{T})$. The theorem also provides specific optimal bounds for $r$
depending on the von Neumann algebra $\mathcal{M}.$ When $\mathcal{M}=M_n(\mathbb{C})$, the optimal bound for $r$ is at most $\frac{n}{3n-2}$ for any general $n$, $1/2$ for $n=2$, and $\sqrt{2}-1$ for $n=3.$ When $\mathcal{M}=B(\mathcal{H})$ for an infinite dimensional Hilbert space $\mathcal{H}$ or a commutative von Neumann algebra, the optimal bound is $1/3$.

Lata and Singh's work also provides versions of Bohr's inequality under more relaxed hypotheses. For the class $\mathcal{C}_1$, they proved that an inequality similar
to (\ref{bohr2}) can be obtained for a larger collection of operators. Additionally, for finite matrices in $M_n(\mathbb{C})$, they showed that the optimal bound for
$r$ is $1/3$ for all $n\ge 2$, unlike the specific bounds given by Theorem \ref{bohr1}, as mentioned above. This was achieved by constructing specific non-trivial
matrices that cause the inequality to fail for $r>1/3$.

\vspace{.3 cm}

\noindent{\bf The Gleason--Kahane--\.Zelazko (GKZ) Theorem}

\vspace{.2 cm}

The GKZ theorem \cite{Gle, KZ, Zel} characterizes multiplicative linear functionals on Banach algebras. Lata and Singh, along with Jaikishan, developed new analogues for function spaces in a vein similar to that pursued by Mashreghi and Ransford \cite{MR1}, but employing new methods in different directions. In \cite{JLS}, they established a general GKZ theorem on the vector space of complex polynomials, showing that multiplicative linear functionals–those with $F(1)=1$ and $F((z-\alpha)^n)\neq 0$ for all $n$ and $\alpha$ outside a fixed open disc–are precisely point evaluations. \emph{Mathematical Reviews} described this as a far-reaching generalization. Their work further extends to vector spaces with topologies not tied to algebraic structure.

\begin{theorem}\label{JLS}\cite[Theorem B]{JLS} Let $\mathcal X\subset \mathcal{F}(r,\mathbb{C})$ be a complex vector space equipped with a topology that satisfies the following properties:
\begin{enumerate}
\item for each $w\in D(0,r)$, the map $f \mapsto f(w)$ is continuous;
\item $\mathcal X$ contains the set of complex polynomials as a dense subset.
\end{enumerate}
Let $\mathcal Y=\{(z-\lambda)^n: n\ge 0, \ \lambda\in \mathbb C, \  |\lambda|\ge r\}$. If $F: \mathcal X \to \mathbb{C}$ is a continuous linear functional such that $F(1)=1$ and $F(g)\ne0$ for all $g\in \mathcal Y$, then there exists $w \in D(0,r)$ such that
\begin{equation*}
F(f)=f(w)
\end{equation*}
for all $f\in \mathcal X$.
\end{theorem}

The authors observe that while their theorem parallels Mashreghi and Ransford’s results, it relies only on a narrow class of outer
functions–powers of linear outer polynomials–and, unlike their proofs, does not use the GKZ theorem.

Lata and Singh, with Jaikishan, extended the Kowalski--$\text{S\l{}odkowski}$ (KS) theorem \cite{KS}, a generalization of the GKZ theorem, to broader topological spaces in \cite{JLS1}. They showed that non-zero continuous functions with values in a Banach algebra satisfying KS-type conditions can be expressed as a composition of a point evaluation and a multiplicative linear functional. Their analogue of the KS theorem stands apart from related works of Sampat \cite{Sam} and Sebastian--Daniel \cite{SD}. Their result provides a several-variables analogue of the KS theorem for vector-valued functions, while Sampat’s yields a GKZ version and Sebastian--Daniel's is confined to Hardy spaces on the unit disc. The proofs given in \cite{JLS1} are elementary, in contrast to the function-theoretic approach of \cite{Sam} and the Hardy space approach of \cite{SD}. They assume the map is non-zero only on certain polynomials, introduce an essential continuity condition, and extend the result beyond analytic function spaces.

The main thrust of our plans for the  future shall centre around–in addition to developing new applications and generalizations of Bohr's inequality as well as the GKZ
theorem–investigations of invariant subspaces of operators on the spaces BMO and VMO and on abstract function spaces associated with compact Abelian groups and group algebras.

\section{von Neumann Algebras}

\noindent \textit{By Kunal Krishna Mukherjee [kunal@iitm.ac.in] from IIT Madras.}

\bigskip

In India, the community of researchers working in von Neumann algebras is led by the author alongside P. Bikram, I. Patri, K. Bakshi, V. Gupta, V. Kodiyalam, S. Ghosh, and their students and postdoctoral collaborators. Broadly speaking, the first three--Mukherjee, Bikram, and Patri—have advanced the analytical, algebraic, and dynamical aspects of von Neumann algebras, while the latter group--Bakshi, Gupta, Kodiyalam, and Ghosh--have primarily concentrated on $\mathrm{II}_1$ factors and subfactor theory. It is important to emphasize that the development of von Neumann algebras in India reflects a clear generational leap: the field has evolved from an early focus almost exclusively on $\mathrm{II}_1$ factors and subfactor theory to a much broader engagement with factors of all types. This shift marks a significant success story in the Indian operator algebra community, positioning it to contribute meaningfully to global research trends. This article, in particular, concentrates on the analytical, algebraic, and dynamical aspects of von Neumann algebras, without restricting attention to any specific type of factor, thereby reflecting this broader and more mature research landscape. Here, we highlight several of the major results obtained over the past decade, which together illustrate both the depth and the breadth of these investigations.

The recent paper \cite{BCCMSW} studied the relationship between the dynamics of the action $\alpha$ of a discrete group $G$ on a von Neumann algebra $M$, and structural properties of the associated crossed product inclusion $L(G) \subseteq M \rtimes_\alpha G$, and its intermediate  subalgebras. The underlying group is only assumed to be {\em discrete}, regardless of its cardinality, and the ambient von Neumann algebra is just assumed to be $\sigma$-finite. A key tool in the setting of a noncommutative dynamical system is the set of {\em quasinormalizers} for an inclusion of von Neumann algebras. It was shown that the von Neumann algebra generated by the quasinormalizers captures analytical properties of the inclusion $L(G) \subseteq M \rtimes_\alpha G$, such as the Haagerup Approximation Property, and is essential to capturing {\em almost periodic} behaviour in the dynamical system. This yields a new description of the Furstenberg-Zimmer distal tower for an ergodic action on a probability space. New versions of the Furstenberg-Zimmer structure theorem for general, tracial $W^*$-dynamical systems were obtained.  This work builds upon more than eighty years of profound developments in the theory of operator algebras, with Dixmier’s seminal contributions in the 1950s (notably Dixmier, 1954) marking the starting point of systematic investigations into the structure of von Neumann algebras. The study of quasinormalizers, in particular, has revealed a remarkably rich and intricate structure that continues to pose challenging open questions. We believe that, despite the substantial progress achieved to date, there remains a long trajectory ahead before a comprehensive understanding of these objects is reached.

A long paper \cite{DM} (in a series of five) takes up a comprehensive study of uniformly left bounded $($resp. left-right bounded$)$ orthonormal bases in GNS spaces of infinite-dimensional von Neumann algebras in the framework of both faithful normal states and $f.n.s.$ weights. Necessary and sufficient conditions on a closed subspace of a GNS space were provided to guarantee the existence of an orthonormal basis of uniformly left bounded $($resp. left-right bounded$)$ vectors. In the context of states, while a basis of the first kind exists for all GNS spaces, $\mathbf{B}(\ell^2)$ (and purely atomic algebras in general) is excluded for a basis of the latter kind. However, in the context of weights, there are no such aforesaid obstructions. In the context of weights, the GNS space of every infinite-dimensional von Neumann algebra admits a uniformly left and right bounded orthonormal basis such that the aforesaid bound is arbitrarily small.

If $M$ is an infinite-dimensional factor and $\varphi$ is a faithful normal state on $M$, then given $\epsilon>0$, the associated GNS space admits a uniformly left bounded orthonormal basis $\mathcal{O}$ such that
$\sup_{\xi\in\mathcal{O}}\|{L_\xi}\|\leqslant (1+\sqrt{2})+\epsilon$. Similar statements were obtained for left-right bounded bases when $M$ is either of type $\mathrm{II}$ or $\mathrm{III}_\lambda$ with $\lambda\in [0,1)$.
As an outcome, it follows that the GNS space of any type $\rm{II}_1$ factor with respect to its tracial state admits an orthonormal basis consisting of images of self-adjoint operators from the ball of radius $(1+\sqrt{2})+\epsilon$, for every $\epsilon>0$. Related questions regarding orthounitary bases remain open and untouched since the Baton-Rouge conference '67.

The factoriality of $q$-deformed Araki-Woods algebras (algebras associated with deformed harmonic oscillator) which was open since 2001 after repeated abortive attempts by stellar researchers in the field was settled except for one case in \cite{BM} and \cite{BMRW}. This led to stimulus for others because recently the last case has been settled by others. The constructions of such algebras were first dreamt of in 1974 to create field theories that allow `small violations' of Pauli's exclusion principle. Since such algebras are of type $\rm{III}$ and the generators generate abelian algebras which usually lack appropriate conditional expectation,  there is little to perform calculations. This difficulty was bypassed by choosing an appropriate generator that generates a maximal abelian subalgebra with conditional expectation and that which interacts with the other generators so as to allow calculation of associated bimodules. This was the first calculation of masa-bimodules in the type $\rm{III}$ setting, the $S$-invariant of Connes was calculated successfully and it was shown that in the case the factor turns out to be type $\rm{III}_1$, they saisfy the bicentralizer conjecture of Connes. Several subsequent papers, both in India and internationally, have pursued the study of deformed Araki–Woods von Neumann algebras. However, this represents a noteworthy instance where the direction of inquiry was not one of following existing trends but, rather, one in which our work set the precedent and was later followed by others.

A systematic and comprehensive study of joinings of non-commutative dynamical systems was carried out in a series of 3 papers \cite{BCM1,BCM2,BCM3}. The study of actions of non-commutative groups on von Neumann algebras was extended to a significant height, strengthening a well-known result of Hoegh-Khron, Landstad and St{\o}rmer on actions of compact groups on von Neumann algebras that appeared in Ann. Math. in '82. The result is as follows: If a locally compact group $G$ acts ergodically on von Neumann algebra $M$ with separable predual preserving a faithful normal state $\varphi$, then any finite-dimensional invariant subspace of the associated Koopman representation in the GNS Hilbert space is contained inside the image of the centralizer $M^\varphi$. Consequently, if $M^\varphi=\mathbb{C}1 $ i.e., $M$ is a $\rm{III}_1$ factor, then any such ergodic action is weak mixing. It basically means non-commutative dynamics are chaotic on large regions of the phase space.

A comprehensive study of compact bicrossed products, a construction in the theory of Quantum Groups, was carried out in \cite{FMP}  and an infinite family of discrete quantum groups with property (T) (of Khazdan and Margulis) were exhibited, a single example of which was unknown then.

The paper in \cite{MP}  initiated the study of non-commutative dynamical systems of the form $(\mathbb{G},\Gamma)$, where a discrete group
acts on a compact quantum group (CQG)
by quantum automorphisms. Combinatorial conditions for such dynamical systems to be ergodic, mixing, compact, etc., were obtained, and a wide variety of examples to illustrate these conditions were provided. A well-known theorem of Halmos to demonstrate ‘reversal of arrows’ in the ergodic hierarchy relevant to the context was generalized. Also, a study of spectral measures for actions of (non-commutative) groups was made. An investigation of the structure of such dynamical systems was made, and under certain restrictions, the existence and uniqueness of the maximal ergodic invariant normal subgroup of such systems was established. As an application, the size of normalizing algebras of masas arising from groups in von Neumann algebraic CQGs was studied, and it was shown that the normalizing algebra of such masas is the von Neumann algebra generated by co-commutative CQGs.

\section{Applications of Birkhoff-James orthogonality in studying the geometry of operators between Banach spaces}

\noindent \textit{By Kallol Paul [kalloldada@gmail.com] from Jadavpur University; and Debmalya Sain [saindebmalya@gmail.com] from IIIT Raichur, Karnataka}

\bigskip

Birkhoff-James orthogonality is arguably the most natural and well-studied concept of orthogonality in the setting of Banach spaces. Introduced by Birkhoff \cite{Birkhoff} and later developed by James \cite{James}, it provides key insights on the structure of Banach spaces and perhaps more importantly, bounded linear operators between them. \emph{For elements $x,y$ in a Banach space $\mathbb{X},$ $x$ is Birkhoff-James orthogonal to $y$} ({\emph{denoted as $x \perp_B y$}) \emph{if $\|x + \lambda y \| \geq \|x\|,$ for all scalars $\lambda.$}  A significant part of  our research in the last two decades is dedicated to exploring bounded linear operators by employing orthogonality techniques. Letters $\mathbb{X},\mathbb{Y}$ denote Banach spaces and $\mathbb{H}$  a Hilbert space. The dual space of $\mathbb{X}$ is denoted by $\mathbb{X}^*$.  $B_\mathbb{X}$ and $S_{\mathbb{X}}$ stand for the unit ball and the unit sphere of $\mathbb{X},$ respectively.  Let $\mathbb{L}(\mathbb{X}, \mathbb{Y})$ $( \text{resp.} ~ \mathbb{K}(\mathbb{X}, \mathbb{Y}))$ be the space of all bounded (resp. compact) linear operators between $\mathbb{X}$ and $\mathbb{Y}.$ The space of all bounded (resp. compact) linear operators  on $\mathbb{X}$ is denoted by $\mathbb{L}(\mathbb{X}) $ (\text{resp.} ~ $\mathbb{K}(\mathbb{X})$). Given $T \in \mathbb{L}(\mathbb{X}, \mathbb{Y})$, the norm attainment set of $ T, $ denoted by $M_T,$ is defined as $M_T= \{ x \in S_{\mathbb{X}}: \|Tx\|=\|T\|\}.$ A useful characterization of Birkhoff-James orthogonality in   $\mathbb{L}(\mathbb{H})$, where $\mathbb{H}$ is an Euclidean space, was obtained in \cite{BhatiaSemrl, Paul}. It is worth mentioning that the following conjecture in \cite{BhatiaSemrl} was one of the key motivators for our study in this direction:\\

``\textit{If $T, A \in \mathbb L (\mathbb X), $ where  $\mathbb X$ is a finite-dimensional Banach space, are such that $ T \perp_B A$ then there exists $ x \in M_T $ such that $ Tx \perp_B Ax.$}''\\

In  this context, our main contribution in \cite{SainPaul} is to illustrate the importance of the norm attainment set $ M_T $ in the whole scheme of things. More precisely, we show that for real Banach spaces, the above conjecture is true whenever $M_T= \pm D,$ where $D$ is a connected set. On the other hand, easy counter-examples to the above conjecture can be constructed by considering $ M_T $ to be not of this form. This insight also allowed us to  characterize Euclidean spaces: ``\emph{A finite-dimensional real Banach space $\mathbb{X}$ is a Hilbert space if and only if for any $T \in \mathbb{L}(\mathbb{X})$, $M_T$ is the unit sphere of some subspace of $\mathbb{X}.$}" Furthermore, in \cite{SainPaulMal18}, we extended the characterization to compact operators   $T \in \mathbb{K}(\mathbb{X}, \mathbb{Y}),$ where $\mathbb{X}$ is a reflexive Banach space. In the most general setting of  $\mathbb{L}(\mathbb{X}, \mathbb{Y}),$ where norm attainment is not a priori guaranteed, complete characterizations of Birkhoff-James orthogonality have been obtained in \cite{MalPaulSain} and \cite{SainPaulMal18}, with important applications to the study of Gateaux differentiability (also called smoothness in the language of geometry) in operator spaces.\\
	
Let us recall that an element $x \in \mathbb{X}$ is said to be smooth if $\{ f \in S_{\mathbb{X}^*}: f(x)=\|x\|\}$ is a singleton. Smoothness of an element is fundamentally related to orthogonality in the following way: A non-zero vector \emph{$x \in \mathbb{X}$  is smooth if and only if Birkhoff-James orthogonality at $x$ is right additive, i.e., $x \perp_B y$ and $ x \perp_B z \implies x \perp_B (y+z).$} By using the previously mentioned characterizations of Birkhoff-James orthogonality, we studied the smoothness of bounded linear operators. In \cite{PaulSainGhosh}, complete characterizations of smooth operators in $\mathbb{L}(\mathbb{H})$ and in $\mathbb{K}(\mathbb{X})$ have been obtained, for a reflexive real Banach space $\mathbb{X}.$ These results were further generalized in \cite{SainPaulMalRay}, where the smooth operators in $\mathbb{L}(\mathbb{X}, \mathbb{Y})$ have been characterized. We have also explored a generalized notion of smoothness, namely, the $k$-smoothness. An element $x \in S_{\mathbb{X}}$ is said to be $k$-smooth if $dim~ span \{ f \in S_{\mathbb{X}^*}: f(x)= \|x\|\}=k.$ A complete identification of $k$-smooth operators in $\mathbb{L}(\mathbb{H})$ has been obtained in \cite{MalDeyPaul}.  Some related results have also been presented for operators on $\ell_1^n$ and $\ell_\infty^n.$ We have further extended this study in \cite{SainSohelPaul}, where we have solved the $k$-smoothness problem for operators between polyhedral Banach spaces. \\
	
We have also explored another fundamental geometric notion in the space of bounded linear operators, namely, the extreme contractions. We recall that extreme contractions are simply the extreme points of the unit ball of the space of all bounded linear operators. Although the extreme contractions in $\mathbb{L}(\mathbb{H})$ have been characterized nearly seven decades ago, a tractable solution in a general Banach space setting remains elusive, even in the finite-dimensional case. In \cite{SainPaulMal21}, we completely characterized these special points in $\mathbb{L}(\mathbb{X}, \mathbb{Y}),$ where $\mathbb{X}, \mathbb{Y}$ are two-dimensional strictly convex and smooth Banach spaces. We further study this problem for finite-dimensional polyhedral spaces in \cite{SainRayPaul}, where it was shown that \emph{if $\mathbb{X}$ is an $n$-dimensional polyhedral Banach space and $T \in \mathbb{L}(\mathbb{X}, \mathbb{Y})$ is an extreme contraction, then $M_T$ contains $n$ linearly independent extreme points of
$B_{\mathbb{X}}.$} Throughout these studies, Birkhoff-James orthogonality has played a central role.\\

More recently, we have applied our work on Birkhoff-James orthogonality to the theory optimization in Banach spaces, focusing on classical best approximation and best coapproximation problems. In \cite{Sain}, best approximations to compact operators between Banach spaces and Hilbert spaces have been studied from the view point of Birkhoff-James orthogonality. As an application, some distance formulae have been presented in the space of compact operators. Later this work has been further extended and generalized in \cite{MalPaul}, where we have obtained several distance formulae for a compact operator to a subspace. The best coapproximation problem, a much less explored but equally natural counterpart of the classical best approximation problem, was studied in \cite{SainSohelGhoshPaul}. We have completely characterized best coapproximation points from an $n \times n$ matrix to a subspace of diagonal matrices and moreover, identified the coproximinal and the co-Chebyshev subspaces in this setting, which are important for many practical purposes.\\
	
In a nutshell, our research on Birkhoff-James orthogonality and its various applications in the isometric theory of Banach spaces illustrate the importance of orthogonality in understanding the analytic as well as the geometric properties of Banach spaces and of bounded linear operators between them. In the last decade, several groups of Banach space theorists have become interested in applying orthogonality techniques in the setting of operator spaces. We strongly expect this trend to continue for good reasons and it is our intention to further explore orthogonality related topics in the more general framework of topological vector spaces.

\section{Commutators and Lie algebras of compact operators}

\noindent \textit{By Sasmita Patnaik [sasmita@iitk.ac.in] from IIT Kanpur.}

[Dedicated to the memory of Gary Weiss]

\bigskip

This article centers around two long-standing open questions spanning over four decades: the Pearcy--Topping's compact commutator problem \cite[Section 2, Problem 1]{Pearcy-Topping} and the Wojty$\acute{\text{n}}$ski's Lie algebra simplicity problem \cite[Question 3]{Wojtynski}. We briefly summarize the progress made with emphasis on the approach to these problems.

Commutators are operators of the form $AB- BA$, where $A$ and $B$ are bounded operators on a complex Hilbert space. When the Hilbert space is finite dimensional, then a matrix is a commutator if and only if its trace is zero \cite{Shoda}. The situation changes in the infinite-dimensional Hilbert space setting. In 1965, Brown--Pearcy characterized all commutators of bounded operators: a bounded operator is a commutator if and only if it is not a nonzero scalar perturbation of a compact operator \cite[Theorem 3]{Brown-Pearcy}. Later, Pearcy--Topping initiated the study of commutators of compact operators by asking four seminal questions on its structure. Among them, the following remains unresolved \cite[Section 2, Problem 1]{Pearcy-Topping} (1971): is every compact operator a commutator of compact operators? In this direction, Anderson constructed infinite block matrices whose commutator is a rank-one projection operator (infinite kernel) \cite[Theorem 1]{Anderson}, thereby affirmatively answering the key test question posed by Pearcy--Topping in \cite{Pearcy-Topping}. This was a groundbreaking work of Anderson and the starting point of our research in the subject. By making a modest modification of Anderson's matrix model, we (joint with Gary Weiss) answered a 35-year-old test question posed by Gary Weiss \cite[Theorem 3.1]{Daniel-Sasmita-Gary}: there is a rich class of strictly positive (zero kernel) compact operators with \textit{repeating} eigenvalues that are commutators of compact operators.

The work mentioned in this paragraph is a joint work with Jireh Loreaux, John Petrovic, and Gary Weiss. In continuation of our earlier work on commutators, by making a controlled perturbation of Anderson's matrices and extending the techniques of Anderson's matrix model, we further enlarged the known class by constructing a family of strictly positive compact operators with \textit{distinct} eigenvalues \cite[Theorem 2.10]{Loreaux-Patnaik-Petrovic-Weiss}. After the successful implementation of certain variations of Anderson matrix model, we encountered obstacles in the use of Anderson's matrix model that had to do with the arithmetic growth size of its block matrices \cite[Theorem 3.4, Corollary 3.5]{Loreaux-Patnaik-Petrovic-Weiss}. A deeper analysis of Anderson's matrix model opened a new vein in our study of commutators - sparsification of matrices in universal block tridiagonal matrix forms  \cite[Theorem 20.4, Theorems 20.7-20.8]{Patnaik-Petrovic-Weiss}, which is a joint work with John Petrovic and Gary Weiss. We note for the reader that any progress that has been made so far on this problem involved alterations in the Anderson's matrix model - a very matricial approach. After getting soaked in the matrix constructions for a while, we next viewed this problem from a different angle independent of matrix representations. We gave a fresh new perspective to the Pearcy-Topping problem by presenting more general constraints on the solutions of the commutator equation $CZ-ZC = T$ in question. They are of two different types; quantitative in terms of $s$-numbers \cite[Proposition 4.1]{Loreaux-Patnaik-Petrovic-Weiss} and qualitative in terms of principal ideals generated by the operators $T,C$, and $Z$ \cite[Corollary 4.8]{Loreaux-Patnaik-Petrovic-Weiss}. Several questions popped up in this direction which one could take up to explore in \cite[Section 4]{Loreaux-Patnaik-Petrovic-Weiss}.

Our work over a decade on commutators of compact operators has led to various offshoots that adds perspective in several directions. For instance, Jireh Loreaux and Gary Weiss found connecting threads between our study of single commutators of compact operators with the concept of generalized traces (beyond spectral traces) on operator ideals \cite[Paragraph following Theorem 9]{Loreaux-Weiss}; and with the splitting of infinite matrix representations of compact operators and the membership of its parts in operator ideals \cite[Paragraph 6]{Loreaux-Weiss}. Our universal block tridiagonal matrix form associated with a compact operator allowed us to reformulate an equivalent Pearcy--Topping question as follows: a compact operator is a commutator of compact operators if and only if one can solve a certain infinite system of finite matrix equations. Solving an infinite system of finite matrices with varying sizes has its own challenges, but this equivalence broadly links the study of elementary operators with the study of finite matrix theory.

The commutator operation $AB - BA$ on the space of bounded linear operators gives a non-associative structure to it, which is called a Lie algebra \cite[Definition 2.5]{Patnaik}. A concept in Lie algebra that inherits both an algebraic and topological property is the topological simplicity of Lie algebras (i.e., if it has no nontrivial closed Lie ideals). The Lie algebra of bounded operators is not topologically simple because the closed Lie ideal of compact operators is a nontrivial Lie ideal. In 1976, Wojty$\acute{\text{n}}$ski proved that if the Banach-Lie algebra of compact quasinilpotent operators contains a nonzero finite-rank operator, then the Lie algebra is not topologically simple \cite[Theorem 6]{Wojtynski}. And in that spirit, he raised the question \cite[Question 3]{Wojtynski}: Does every Banach-Lie algebra of compact quasinilpotent operators on a Banach space contain a nontrivial closed Lie ideal?  Some partial results centered around this question were obtained by Bre$\check{\text{s}}$ar, Shulman, and Turovskii in \cite{Bresar-Shulman-Turovski}. To the best of our knowledge, in any of the developments made so far on this problem, the consideration of the adjoint representation of Lie algebras in the study of the simplicity of Lie algebras seemed inevitable.

In 2022, the author of this article offered a new and direct hands-on approach to the study of simplicity of Lie algebras of compact operators acting on an infinite-dimensional  separable complex Hilbert space that avoids the adjoint representation technique. We achieved this goal by first focussing on the study of the algebraic simplicity of Lie algebras of compact operators by addressing the question  \cite[Section 1, Paragraph 2]{Patnaik}:  Which infinite-dimensional Lie algebras of compact operators on a Hilbert space are algebraically simple?
 By employing the notion of soft-edged ideals in the theory of operator ideals, we obtained a sufficient condition that guarantees the non-simplicity of Lie algebras of compact operators \cite[Theorem 3.5 and Corollary 3.8]{Patnaik}.
Soft-edged ideals were first introduced by Kaftal and Weiss and these ideals have played a significant role in various aspects of operator ideal theory, see \cite{Kaftal-Weiss2007}-\cite{Kaftal-Weiss2002} and the references therein. This approach establishes a bridge connecting a purely algebraic problem (simplicity) with an analytic behaviour of certain operator ideals (ideal-softness). We also emphasized that the investigation of simplicity of Lie algebras of finite-rank operators requires a separate treatment \cite[Question]{Patnaik}. Ultimately, we intend to address the Wojty$\acute{\text{n}}$ski question in the Hilbert space framework hoping that the operator ideal techniques might shed some light on the problem unlike the intractability of the problem in more general Banach spaces.

\section{On subclasses of norm attaining operators on Hilbert spaces}

\noindent \textit{By G. Ramesh [rameshg@math.iith.ac.in] from IIT Hyderabad}

\bigskip

Norm attaining operators on Banach spaces have been studied extensively for several decades. When considered on Hilbert spaces, these operators exhibit richer structural and geometric properties. In recent years, significant attention has been devoted to both norm attaining operators and their subclasses, as well as to minimum attaining operators and their subclasses. In particular, both classes form dense subsets of $\mathcal B(H)$, the Banach space of all bounded linear operators on a Hilbert space $H$, equipped with the operator norm.

In this note, we focus on two central themes in operator theory: the \emph{norm attaining property} and the \emph{invariant subspace property} for certain subclasses of norm attaining and minimum attaining operators. We describe these details in the following.

Recall that a closed subspace $ M$ of $ H$ is said to be an \textbf{invariant subspace} for $ T$ or \textbf{invariant} under $ T$ if $ T(M) \subseteq M$; that is, $ Tx \in M$, whenever $ x \in M$. Moreover, $ M$ is said to be \textbf{reducing} for $ T$ if both $ M$ and its orthogonal complement $ M^{\perp}$ in $ H$ are invariant under $ T$. If $ M$ is invariant for every $ S \in \mathcal{B}(H)$ that commutes with $ T$, then $ M$ is said to be \textbf{hyperinvariant} under $ T$. The invariant subspace (respectively, hyperinvariant subspace) problem can be stated as follows:
\begin{problem}
Let
$H$ be a separable, infinite-dimensional Hilbert space. Does every
$T\in \mathcal B(H)\setminus \mathbb{C}I$ admit a nontrivial invariant (respectively, hyperinvariant) subspace?
\end{problem}

\smallskip

\noindent \textsf{Absolutely Norm attaining Operators:} We say  $T\in \mathcal B(H)$ is norm attaining if there exists $x\in H$ with $\|x\|=1$ such that $\|Tx\|=\|T\|$ and is absolutely norm attaining (or $\mathcal{AN}$-operator, for short) if for every non-zero closed subspace $M$ of $H$, the restriction operator $T|_{M}:M\rightarrow H$ is norm attaining. We denote the set of all absolutely norm attaining operators on $H$ by $\mathcal{AN}(H)$.

The class of \emph{absolutely norm attaining operators} was introduced in \cite{carvajalneves1} based on the properties of compact operators and isometries, with the motivation of studying the invariant subspace problem. It includes compact operators, isometries, and partial isometries with finite-dimensional kernels.

Several characterizations of positive $\mathcal{AN}$-operators are provided in \cite{PP,GRpara,GRSSS1,GRVN}, while characterizations for normal $\mathcal{AN}$-operators appeared in \cite{NBGRBulsci,GRpara,GRVN}. The self-adjoint case is treated in \cite{NBGRBulsci,GRVN}.

A central motivation for studying this class is its connection to the invariant subspace problem, which can be approached via various representations of such operators. For instance, the spectral decomposition of absolutely norm attaining normal operators is discussed in \cite{NBGRBulsci,GRpara}, and the hyponormal case is considered in \cite{NBGRBulsci}. A key question here is: \emph{When does a hyponormal operator become normal?} It is well known that compact hyponormal operators are normal. In \cite{NBGRBulsci}, it is shown that a hyponormal $\mathcal{AN}$-operator whose Weyl spectrum coincides with its essential spectrum must be normal. This result has been generalized to the operator norm closure of $\mathcal{AN}$-operators in \cite{GRSSS4}.  For developments on the functional calculus of $\mathcal{AN}$-operators, we refer to \cite{GHUY1}, and for results on the stability of $\mathcal{AN}$-operators under the induced Aluthge transform, we refer to \cite{GHUY2}.

Another noteworthy subclass is the set of bounded operators attaining their norm on every non-zero reducing subspace. This class does not contain the operator norm closure of the $\mathcal{AN}$-class, nor does it contained the operator norm closure of $\mathcal{AN}(H)$. However, it is dense in $\mathcal B(H)$ with respect to the operator norm. In particular, every bounded operator can be approximated by an operator in this class with a nontrivial invariant subspace (see \cite{GROHSSJOT} for details).

\smallskip

\noindent \textsf{Absolutely Minimum attaining Operators:} A natural counterpart to the class of  $\mathcal{AN}$-operators is the class of \emph{absolutely minimum attaining operators}. Recall that $T \in \mathcal B(H)$ is \emph{minimum attaining} if there exists $x \in H$ with $\|x\|=1$ such that $\|Tx\| = m(T)$, the \emph{minimum modulus} of $T$, where
\[
m(T) = \inf \{\|Tx\| : x \in H, \; \|x\| = 1\}.
\]
It has been shown that the set of minimum attaining operators is dense in $\mathcal B(H)$ with respect to the operator norm \cite{SHKGRMNM}, and moreover, a Bishop--Phelps--Bollobás type theorem holds for this class \cite{NBGRBPBP}.  We say $T\in \mathcal B(H)$ to be \textit{absolutely minimum attaining} (or $\mathcal{AM}$, for short) if for every non-zero closed subspace $M$ of $H$, the restriction operator $T|_{M}:M\rightarrow H$ is minimum attaining. The collection of all such operators is denoted by $\mathcal{AM}(H)$, a class first introduced in \cite{carvajaneves2}.

The behaviour of $\mathcal{AM}$-operators contrasts sharply with that of $\mathcal{AN}$-operators. For instance, every $\mathcal{AM}$-operator has closed range, whereas an infinite-rank compact operator can be $\mathcal{AN}$ without having closed range. Characterizations and properties of positive $\mathcal{AM}$-operators are given in \cite{GRS18}, and several spectral characterizations appear in \cite{NBGRspectral}.

The existence of hyperinvariant subspaces for Toeplitz and Hankel operators in both $\mathcal{AN}$ and $\mathcal{AM}$ classes has been studied in \cite{GRSSS2, GRSSS3}. Interestingly, the operator norm closures of $\mathcal{AN}(H)$ and $\mathcal{AM}(H)$ coincide. Specifically,
\begin{align*}
\overline{\mathcal{AN}(H)}& = \overline{\mathcal{AM}(H)} \\
                        &= \{K + \alpha V: K \; \text{compact}, \; V \; \text{partial isometry with }, \; \dim (\ker(V))<\infty, \; \alpha \ge 0\}.
\end{align*}
For further details, see \cite{GRSSS1}. This class includes compact perturbations of unitaries, an important class of operators in operator theory. The existence of invariant (hyperinvariant) subspaces for this class remains an open problem.  In this direction, recent progress \cite{NBGRHISP} shows that if $T \in \mathcal B(H)$ satisfies $\sigma_{\text{ess}}(T^*T)=\{\alpha\}$ for some $\alpha\geq 0$ and the sequence $(\alpha-\lambda_n) \in \ell^p(\mathbb N)$ for some $1 < p < \infty$, where $\lambda_n\in \sigma_d(T^*T)$,  then $T$ has a nontrivial hyperinvariant subspace. Here, $\sigma_{\text{ess}}(A)$ and $\sigma_{d}(A)$ denote the essential and discrete spectra, respectively, of a self-adjoint operator $A\in \mathcal B(H)$. As a consequence we can easily recover a well known result: If $T \in \mathcal B(H)$ with $I - T^*T\in \mathcal K_p(H)$ ( the Schatten-$p$ class) for some $1 < p < \infty$, then $T$ admits a non-trivial hyperinvariant subspace.

For a comprehensive survey on the spectral theory of unbounded absolutely minimum-attaining operators, we refer the reader to \cite{SHKGRsurvey}.

The first step towards generalizing $\mathcal{AN}$-operators to Banach space operators is to obtain a representation of compact operators on Banach spaces, as established in \cite{RVS}.

\section{Generalized Parrott homomorphisms and Banach space geometry}

\noindent \textit{By Samya Kumar Ray [samyaray7777@gmail.com] from The Institute of Mathematical Sciences (HBNI), Chennai; Rajeev Gupta [rajeev@iitgoa.ac.in] from IIT Goa; and Arpita Mal [arpitamalju@gmail.com] from Dhirubhai Ambani University.}

\bigskip

For $n \in \mathbb{N}$, let $\mathcal{C}_n$ be the set of $n$-tuples of commuting contractions and $\mathbb{C}_k[z]$ be the set of polynomials of degree at most $k$. Define
\[
\mathcal{C}_k(n) := \sup \big\{ \|p(T)\| : \|p\|_{\infty,\mathbb{D}^n} \leqslant 1,\, p \in \mathbb{C}_k[z],\, T \in \mathcal{C}_n \big\},
\qquad
\mathcal{C}(n) := \lim_{k \to \infty} \mathcal{C}_k(n).
\]
We have $\mathcal{C}(1)=\mathcal{C}(2)=1$ by dilation theorem of Sz.-Nagy \cite{SzNagy2010} and And\"{o} \cite{Ando} respectively. For $n \geqslant 3$, $\mathcal{C}(n)$ grows faster than any power of $n$ \cite{Gacre18}, but it is unknown whether $\mathcal{C}(3)$ (or any $\mathcal{C}(n)$, $n \geqslant 3$) is finite \cite{Pisier2001}. Varopoulos–Kaijser’s \cite{Varopoulos1974} construction shows $\mathcal{C}_2(3)>1$, raising the problem of whether $(\mathcal{C}_2(n))_{n \in \mathbb{N}}$ is bounded.
Varopoulos \cite{Varopoulos1976} proved that
$
K_G^{\mathbb{C}} \leqslant \lim_{n \to \infty} \mathcal{C}_2(n) \leqslant 2K_G^{\mathbb{C}},
$
where $K_G^{\mathbb{C}}$ is the complex Grothendieck constant. Since $K_G^{\mathbb{C}} > 1$, this \emph{Varopoulos inequality} shows the eventual failure of von Neumann’s inequality in high dimensions. He asked whether the limit equals $K_G^{\mathbb{C}}$? In \cite{GuRa18} this question was answered in the negative by proving
$
\lim_{n \to \infty} \mathcal{C}_2(n) \geqslant 1.118\, K_G^{\mathbb{C}}.
$
Arveson \cite{Arveson1998} provided a remarkable equivalent criterion for $\partial \mathbb{D}^3$-normal dilation in terms of contractivity vs. complete contractivity.
Parrott’s example \cite{Parrott1970} showed that a unital contractive homomorphism on $\mathcal{O}(\mathbb{D}^3)$ need not admit a $\partial \mathbb{D}^3$-normal dilation, revealing that complete contractivity can fail in higher dimensions. Here, $\mathcal{O}(\mathbb{D}^n)$ is
the closure in $C(\overline{\mathbb{D}^n})$ of all complex polynomials in $n$ variables. Misra and his coauthors \cite{MisraSastry1984, MisraSastry1986, MisraSastry1987, BBM, Mi94} further studied these examples in detail and extended this to Euclidean balls in $\mathbb{C}^n$ for all $n \geqslant 2$.  Based on these examples, it is natural to study the \emph{Generalized Parrott homomorphisms} described below.

Fix $\ell, m\in \mathbb{N}.$ We denote $(X)_1,$ as the open unit ball of a $d$-dimensional complex Banach space $X.$ Let $J_\ell$ denote the standard Jordan block of size $\ell \in \mathbb{N}$.
and $\mathcal C(\ell,m)$ denote the set of all $d$-tuples of square matrices of size $m (\ell +1)$  of the form $\boldsymbol{T_V^{}}:=(V_1\otimes J_{\ell +1}, \ldots, (V_d\otimes J_{\ell +1}),$ where $\mathbf{V}:=(V_1,\ldots, V_d)$ is a $d$-tuple of commuting matrices of size $m.$ Let us consider the algebra homomorphism $: \Phi_{\ell,\mathbf{V}}:f\mapsto f(\boldsymbol{T_V^{}})$ from $\mathcal{O}((X)_1)$ to $M_{m(l+1)}.$ Define
\begin{equation}
    \alpha(\ell,X):=\sup\{\|\Phi_{\ell,\mathbf{V}}\|_{cb}:\|\Phi_{\ell,\mathbf{V}}\|\leq 1, \mathbf{V}\in \mathcal{C}(\ell,m), m\in\mathbb{N}\}.
    \end{equation}

\begin{Question}
For a fixed $\ell\geq 2,$ when is $\alpha(\ell,X)=1$?
\end{Question}

Note that $\alpha(1,X)$ is nothing but the constant $\alpha(X)$ studied earlier by Paulsen \cite{Pau92}. He observed that if every contractive homomorphism induced by   $\boldsymbol{T_V^{(0)}}$ is completely contractive, then $\alpha(X)=1$ and vice versa.
It is known that $\alpha(\ell_\infty(2)) = \alpha(\ell_1(2)) = 1$, while $\alpha(X) > 1$ whenever $\dim(X) \geqslant 3$, see \cite{Pisier2003, Pau92}.
Moreover, $\alpha(X)=1$ if and only if $X$ admits unique operator space structure \cite{Pau92}. These results highlight the gap between contractivity and complete contractivity, linking dilation theory to the geometry of domains and operator spaces. These connections were further revealed in \cite{BaMi95} where the contractivity vs. complete contractivity of a certain class of Parrott-like homomorphisms (called little Parrott homomorphisms \cite{GuMiRa25}) were shown to be equivalent to the so-called $2$-summing property or Property P. Let us elaborate further. We use the natural notion of positivity for elements of $E \otimes E$, namely, an element
$A \in E \otimes E$ is positive (written $A \geqslant 0$) if it is in the convex hull of the set
of symmetric tensors $x \otimes x$, $x \in E$. Given a finite-dimensional complex Banach space $E$, let
\[
\gamma(E) := \sup \Big\{
   \langle A, B \rangle_{\mathrm{HS}}:  \|A\|_{E \to E^*}\leqslant 1, \, \|B\|_{E^*\to E}\leqslant 1, \, A \geqslant 0, \, B \geqslant 0
\Big\},
\]
where $\langle \cdot, \cdot \rangle_{\mathrm{HS}}$ denotes the Hilbert--Schmidt inner product. Then, in the terminology of \cite{BaMi95}, $E$ is said to possess \emph{Property P} if $\gamma(E) \leqslant 1$. Real Banach spaces with the $2$-summing property have been completely characterized by \cite{ArFiJoSc95}. In \cite{MiPaVaG} the authors proved the remarkable result that a two-dimensional subspace of $(M_2(\mathbb{C}),\|\cdot\|_{op})$ has Property P if and only if it is isometric to $\ell_\infty(2).$ Thereafter, it is natural to ask whether this kind of example covers other familiar complex Banach spaces such as $\ell_p(2)$ for $1\leqslant p\neq 2<\infty$ for which Property P has been well-studied \cite{BaMi95}. Using unitary dilation of contractions it was shown by Gupta and Reza \cite{GuRe18} that $\ell_1(n)$ does not embed into $(M_k(\mathbb{C}),\|\cdot\|_{op})$ for all $n\geqslant 2.$ However, Pisier (see \cite{GuRe18}) raised the question if this result can be generalized to compact operators and G. Misra asked whether $\ell_p(n)$ embeds isometrically into $M_k(\mathbb{C}).$ All these questions were resolved in a series of papers \cite{Ra20}, \cite{ChHoPaPrRa22}, \cite{ChHoPrRa24}; along the way these papers established noncommutative analogues of some well-known results of isometric embeddability between finite dimensional $\ell_p$-spaces. Problems regarding isometric embedding of Schatten-$p$ classes or in general noncommutative $L_p$-spaces are very useful (see \cite{He24, He25}, \cite{Zh20}). The tools developed in \cite{Ra20}, \cite{ChHoPaPrRa22}, \cite{ChHoPrRa24} perhaps can be used to tackle some problems mentioned in \cite{JuPa10}. Moreover, in \cite{Ra20} many two-dimensional Banach spaces were produced without Property P which were outside the scope of \cite{MiPaVaG} and \cite{BaMi95}. However, Despite all this work, it is still not known whether there is a two-dimensional complex Banach space $X$ having Property P but $\alpha(X)>1.$ In $\mathbb{C}^3$, let
$x = (1,0,a)$ and $y = (0,1,b),$
with $|a+b| \vee |a-b| < 1$ but $|a| + |b| > 1$, and set $X = \operatorname{span}\{x,y\}.$
It is known \cite{ArFiJoSc95} that $X$
possesses the $2$-summing property. It would be interesting to know whether $\alpha(X)=1$? This would answer a long-standing open problem \cite{Pisier2003}.

If a space $E$ does not have the $2$-summing property, then it is natural to find the exact value of $\gamma(E).$ Since $\{A\in E\otimes E: \|A\|_{E \to E^*}\leqslant 1,A\geqslant 0\}$ and $\{B\in E\otimes E: \|B\|_{E^* \to E}\leqslant 1,B\geqslant 0\}$ are convex sets, the supremum in $\gamma(E)$ is attained at the extreme points of these sets. Using these it is proved that $\gamma(\ell_1^{\mathbb{R}}(3))=1.125$ \cite{PreGMMR}.
One more interesting observation which goes back to \cite{BaMi95} is that $\sup\{\gamma(\ell_1(n)): n\in\mathbb{N}\}$ is nothing but the positive Grothendieck constant. Motivated by this and some recent work as in \cite{AuPaSzWi20}, \cite{Samthesis} and in \cite{GuMiRa25} we introduced the Grothendieck constant for a pair of Banach spaces.
We prove that for a Banach space $F$ and a sequence of Banach spaces $(E_m)_{m\in\mathbb N}$ with $\dim E_n=n,$ if there is a $C >0$ such that
\[\label{star} \|A\otimes \operatorname{id}_{F}\|_{E_m\check{\otimes}F\to E_n^*\hat{\otimes}F}\leqslant C \|A\|_{E_m\to E_n^*} \tag{$*$}\]  for all $m,n\in\mathbb{N},$ then both $F$ and $F^*$ must have finite cotype. Moreover, assuming that $F$ has the bounded approximation property and that the conjecture in \cite{PisierDuality} has an affirmative answer, we show that $(E_n^*)_{n\geqslant 1}$ satisfies G.T. uniformly \cite{GuMiRa25}. For contractive \emph{little} Parrott homomorphisms $\Phi_{\ell,\mathbf{V}} : \mathcal{O}(\Omega) \to M_{n}$, where $\Omega$ is the dual unit ball of a finite dimensional Banach space $(E,\|\cdot\|)$, we prove the sharp estimate $
\|\Phi_{\ell,\mathbf{V}}\|_{\mathrm{cb}}\leqslant\sqrt{\gamma(E)}.$
This yields a new proof of \cite[Theorem 2.1]{Davidchoi} using the lower bound $K_G^+(\ell_\infty(4),\ell_2(2)) \geqslant 1.1658$, for the complex positive Grothendieck constant for $4\times 4$ positive matrices, obtained in the same paper. It would be interesting to classify Banach spaces satisfying \eqref{star}.

Fix $d,n\in\mathbb N$, $n\geqslant 2.$ Next, for any row vectors $x_1,\ldots, x_d\in\mathbb C^n,$ 
define the operators
\begin{eqnarray}\label{Varopoulos Operators - First kind}
    V_1=\left(
\begin{array}{cccc}
0 & x_1 & 0 & 0\\
0 & 0 & D_{x_1}\otimes I_{n-2} & 0\\
0 & 0 & 0 & x_1^t\\
0 & 0 & 0 & 0
\end{array}
\right), \ldots,
V_d=\left(
\begin{array}{cccc}
0 & x_d & 0 & 0\\
0 & 0 & D_{x_d}\otimes I_{n-2} & 0\\
0 & 0 & 0 & x_d^t\\
0 & 0 & 0 & 0
\end{array}
\right),
\end{eqnarray}
where for each $j=1,\ldots, d,$ $D_{x_j}$ denotes the diagonal operator with diagonal $x_j.$
Note that $\mathbf{V}:=(V_1,\ldots, V_d)$ is a commuting tuple and $V_j$ is a contraction if and only if $\|x_j\|_2\leq 1.$ Let $\mathcal V_{d,n}$ denote the collection of all the tuples $\mathbf{V}$ as defined by \eqref{Varopoulos Operators - First kind} such that $V_j$ is a contraction for each $j=1,\ldots, d.$ Let $C^\sharp(d,n)$ denote the smallest constant $C$ such that \[\|p(\mathbf{V^{}})\|=\Big |\sum_{|\alpha|=n} a_{\alpha}[x_{\alpha_1},\ldots,x_{\alpha_n}]\Big|\leqslant C \|p\|_{\infty, \mathbb D^d},\]
where $[x_{\alpha_1},\ldots,x_{\alpha_n}]=\sum_{j=1}^{n}x_{\alpha_1}^{(j)}\cdots x_{\alpha_n}^{(j)}$, $|\alpha|:=\sum_{j=1}^d\alpha_j$,
for all $\mathbf{V}\in \mathcal V_{\mathbb{D}^d,n}$ and any polynomial $p(z_1,\ldots, z_d)=\sum_{|\alpha|=n}a_\alpha z^\alpha.$ The constant $C^\sharp(d,n)$ is closely related to the higher dimensional Grothendieck constant (cf. \cite{Tonge}).
\begin{Question}
    What is $\lim_{d, n\to \infty}C^\sharp(d,n)$?
\end{Question}
By \cite{Hz-2025} $C^\sharp(d,n)\lesssim_d(\log(n+1))^{d-3}.$

\section{Large time behaviour of  operator semigroups}

\noindent \textit{By Sachi Srivastava [ssrivastava@maths.du.ac.in] from University of Delhi.}
\bigskip

A strongly continuous operator semigroup or $C_0$-semigroup  on  a Banach space $X$  is  a one parameter family $ \{ T_{t}\}_{t \geq 0}$  of bounded linear operators on $X$ satisfying   (i)  $ T_{0} = I_X,$  (ii) $ T_{t} T_{s} = T_{ t + s }, t, s \geq 0 $ (iii) $ t \mapsto T_t x $ is continuous for each $ x \in X. $  Strongly continuous semigroups are important for solving  the abstract Cauchy problem  of the form  $ u'(t) = Au(t), u(0) = x_0 $ where $ A $ is a closed and densely defined linear  operator on a Banach space  $X$.  In fact, this problem is well posed in the sense that for every $x_0 \in X $ it admits a mild solution (given by $u(t)  = T_t x_0 $ ) if and only if $ A$ is the generator of a $C_0$-semigroup  $ \{ T_{ t} \}_{t \geq 0},$ that is,   $ D(A)=\{ x \in X : \lim_{t \downarrow 0} \frac{ T_t (x) - x }{t}  \mbox{ exists in } X \} $ and $ A x = \lim_{t \to 0} \frac{ T_t (x) - x }{t} , x \in D(A) . $ However,  the difficulty  is that in most applications,  the operator $A$ is available, but the semigroup is  not known explicitly.

Driven by applications in evolution equations, control theory and mathematical physics, the study of asymptotic behaviour of $C_0$-semigroups has garnered increasing attention in the last few decades, particularly in the context of stability, non-uniform stability and decay rates \cite{ BOT, BD Sachi, BCT2016, CT2007,  CST20, CPSST23}.  A semigroup    $ \{ T_{t}\}_{t \geq 0}$ is \textit{ uniformly exponentially stable}  if $\| T_t\| \leq M e^{-wt} $  for all $ t > 0 $ where $ M \geq 1, w > 0 $ (that is, the semigroup decays to zero at an exponential rate),  \textit{  strongly stable} if $\| T_t x \| \to 0 $ as $ t \to \infty $ for each $x \in X$ and   if   $ \{ T_{t}\}_{t \geq 0}$ is uniformly bounded it is said to be  \textit{polynomially stable (of order $\alpha)$}  if    there is a constant $ M > 0$ such that $ \| T_t A^{-1}  \|  \leq \frac{M}{ t^{\frac{1}{\alpha}}}  \, \forall \,t > 0.$ Building on work done in \cite{BD Sachi}, Borichev and Tomilov  \cite{BOT} characterised  polynomial decay of the semigroup  via resolvent estimates  for the Hilbert space case. These remarkable results have far reaching ramifications, particularly in the study of  rates of  decay of solutions of  PDE's,  like for example,  the damped wave equation.  While exponential stability behaves as expected under perturbations of the generator of the semigroup,  perturbation results in the context of strong and  polynomial stability have been few and far between except for those appearing in Paunnonen's work \cite{LP13, LP14}.  Interestingly, Rastogi and Srivastava  proved  the following  perturbation result in \cite{RaSr21} :

If    $(B, D(B))$ is the generator of a polynomially stable $ C_{0} $-semigroup $ \{ T_t \}_{t \geq 0}$  with order $ \alpha > 0 $ on a Hilbert space $H$ and $ C $ is a closed linear operator on $H$  such that  $D(B) \subset D(C)$, then $(B+C, D(B))$ generates  a polynomially stable perturbed $C_0$-semigroup $\{T_{B+C}(t)\}_{t\geq 0}$ with order $ \alpha > 0 $ on the Hilbert space $ H $
provided  there exists $ M_{0} \in [0,1) $ such that  (a) $ \sup_{\operatorname{Re} \lambda \geq 0}\left\| CR(\lambda,B) \right\| \leq M_{0} $,  and (b) $ \sup_{\operatorname{Re} \lambda \geq 0}\left\| R(\lambda,B)Cy \right\| \leq M_{0} \left\| y \right\|, \ \forall\, y \in D(C) .$

Furthermore, they applied this and similar results  to conduct an  in-depth analysis of strong and polynomial stability of a particular class of semigroups - the class of delay semigroups. A delay semigroup is a $C_0$-semigroup generated by an opertor matrix of the form  $\mathcal{A} =   \begin{pmatrix}
                  B & \Phi \\
                  0 & \frac{d}{d\sigma}
                  \end{pmatrix} $
associated to an abstract differential equation of the form \begin{center}
\begin{align*}
 \begin{cases}
u^{\prime}(t)  = Bu(t) + \Phi u_{t}, & t > 0,\\
u(0) = x,\\
u_{0} = f, \end{cases}
\end{align*}
 \end{center}
where,  $ \Phi : W^{1,p}([-1,0], X) \rightarrow X $ is a bounded linear operator, called the delay operator $ A : D(A)\subseteq X \rightarrow X $ is the generator of a $ C_{0} $-semigroup $ \{T_t\}_{t \geq 0} $ on a Banach space $X$; $ f \in L^{p}([-1,0], X) $, $ 1 \leq p < \infty $  and
$ u : [-1,\infty) \rightarrow X $, and $ u_{t}:[-1,0] \rightarrow X $ is defined by $ u_{t}(\sigma): = u(t+\sigma) $, $ \sigma \in [-1,0]$, $ t \geq 0 $.  Such equations arise in the study of systems where the current state depends on the past state, for example,  population dynamics and control systems. Rastogi and Srivastava also study the  non-analytic growth bound (first introduced and studied in  \cite{BaBaSr2003}) of the delay semigroup in \cite{RaSr2019}.
The delay operator considered throughout in the above instances is linear and bounded.  A comprehensive framework to address such problems using semigroup theory,  when the delay term in not necessarily linear and nor is the underlying semigroup,  has not been developed yet. This is a part of ongoing work.

An important but basic question in the study of asymptotics of $C_0$-semigroups is whether or not  the exponential growth bound $w_0 (T) $ of the semigroup $ \{ T_t \}_{t \geq 0}$  coincides with the spectral bound $s(A)$ of the generator $A$. This equality is often called the Lyapunov property. In \cite{wei2}  Weis  proved that if  $X =  L^{p}(\mu)$,  where $( \Omega, \mu) $ is a measure space and $ 1< p < \infty $ and $ T_t \geq 0  \forall \,  \, t \geq 0 $  then indeed $ w_0 (T) = s(A).$  In  \cite{Vo2022}, Vogt  gave an elementry proof of this result and along the way proved that in Weis's result positivity of the semigroup can be replaced by (uniform) eventual positivity. Prajapati, Sinha and Srivastava studied this problem for $C_0$-semigroups defined on non-commutative $L^{p}$ spaces for  $ 1 \leq p < \infty $ in \cite{PrSiSr2019}. They gave direct proofs that the property holds for the case $ p =1, 2 $ and also for consistent families of $C_0$-semigroups under added conditions. A complete analog of Weis's result for the non-commutative case remains an open problem.
Further in the context of  non-commutative spaces,  non-commutative versions of  characterisations in terms of generating forms,   for domination of one $C_0$-semigroup by another, where both semigroups act on a classical $L^{2}$ space,
due to Ouhabaaz \cite{Ou1996} and  Barthélemy \cite{Ba1996},  have been obtained by Arora, Chill and Srivastava in \cite{ArChSr2024} for $C_0$-semigroups acting on non-commutative $L^{2}$-spaces.  Yet another strand connecting asymptotics and semigroups defined on non-commutative spaces has been the study of stability of Quantum Dynamical Semigroups undertaken by Bhat and Srivastava \cite{BhSr2015} (for bounded generators) and Kumar, Sinha and Srivastava \cite{KuSiSr2022} for the general setting.  Here,  a one parameter family $  \mathcal{T}= \left \{ \mathcal{T}_t  \right \}_{t \geq 0}$ is a Quantum Dynamical Semigroup  or QDS on $\mathcal{B}(H)$ if  (a) for each $ t \geq 0, \mathcal{T}_{t} : \mathcal{B}(H) \to \mathcal{B}(H) $ is a   completely positive, sub-Markovian, normal map; (b) $\mathcal{T}_{0}(x) = x, \mathcal{T}_{s+t}= \mathcal{T}_{s}\mathcal{T}_{t} $ for all $ s,t \geq 0, x \in \mathcal{B}(H)$ and (c) for each $ x \in \mathcal{B}(H), $  the map $ t \to \mathcal{T}_{t}(x) $ is continuous with respect to the $ weak^{\ast}$ topology of $\mathcal{B}(H).$

Asymptotic anlaysis of positive $C_0$ -semigroups, (that is, $T_t$ maps the positive cone of the underlying space into the positive cone)  is by now  an almost classical area in semigroup theory. In recent years, the area has witnessed a renaissance of sorts with the work of Daners, Gluck and Kennedy  who introduced and  developed an abstract theory of eventually positive $C_0$-semigroups on Banach lattices in a series of papers  starting with \cite{DGK2016,DGK12016,DG2017}.   A $C_0$-semigroup is said to be (uniformly)\textit{ eventually positive} if there exists $t_0 > 0$ such that $ T_t  \geq 0 $ for all $ t \geq  t_0.$  Related notions include  individual eventual positivity, asymptotic positivity  and local eventual positivity (see \cite{Ar22}). In \cite{DG2018}, Daners and Gl\"uck studied the eventual positivity of semigroups under bounded perturbations of the generator and concretely demonstrated that the perturbation theory for eventually positive semigroups is quite different from that of positive semigroups. The case of unbounded perturbations of the generator in this context is treated by  Pappu, Rastogi and Srivastava in \cite{PaRaSr2025}.  The authors then use their results to deduce eventual positivity of delay semigroups under some conditions. Eventually positive  $C_0$-semigroups are a rich and active  subject of research with several open problems.

\section{Recent developments in $E_0$-semigroups and product systems}

\noindent \textit{By S. Sundar [ssundar@imsc.res.in, sundarsobers@gmail.com] from The Institute of Mathematical Sciences (HBNI), Chennai.}

\bigskip

Broadly speaking, $E_0$-semigroups and product systems fall under the broad framework of non-commutative dynamical systems.  The theory of $E_0$-semigroups which is the study of semigroups of endomorphisms on von Neumann algebras was introduced by Powers (\cite{Powers_Index}).    Arveson in his seminal papers (\cite{Arv_Fock}, \cite{Arv_Fock2}, \cite{Arv_Fock3}, \cite{Arv_Fock4}) gave an alternative way of viewing them by introducing product systems. Most of the  papers that appeared during the 80's, 90's and the early 2000's focussed solely (rightfully so) on the $1$-parameter case.  Both the notion of $E_0$-semigroups and product systems  generalise easily to more general semigroups, and it is worthwhile to investigate them  in more generality.

We give a  brief summary of the results obtained in this regard.  Most of the work was done  by the author in collaboration with mathematicians (Anbu Arjunan, S.P. Murugan, R. Srinivasan, Piyasa Sarkar, and C.H. Namitha) who worked/are working in Chennai (CMI/IMSc) during the last ten years.

\begin{center}
\textsf{Main results}
\end{center}

Let $P$ be a measurable semigroup, and let $H$ be a separable Hilbert space. A semigroup $\alpha=\{\alpha_x\}_{x \in P}$ of unital, normal $^*$-endomorphisms of $B(H)$ is called an  $E_0$-semigroup on $B(H)$.  Two $E_0$-semigroups $\alpha$ and $\beta$ on $B(H)$ are said to be cocycle conjugate if $\alpha_x(.)=U_x\beta_x(\cdot)U_x^*$ for some family of unitaries $\{U_x\}_{x \in P}$ such that $U_x\alpha_x(U_y)=U_{xy}$ for every $x,y \in P$. A product system over $P$ is a semigroup $E$ together with a surjective semigroup homomorphism $p:E \to P$ such that the fibres $E(x):=p^{-1}(x)$ are separable Hilbert spaces, and  the map $E(x) \otimes E(y) \ni u \otimes v \to uv$ is a unitary operator for every $x,y \in P$.  The  measurable structures involved are ignored for  this report.

 Given an $E_0$-semigroup $\alpha$, we can attach a family of Hilbert spaces $E_\alpha:=\{E(x)\}_{x \in P}$ as follows: set
\[
E(x):=\{T \in B(H): \alpha_x(A)T=TA~~ \forall A \in B(H)\}.\]
Then, $E_\alpha:=\{E(x)\}_{x \in P}$ is a product system over $P$ with the product rule  given by composition, and with the inner product on $E(x)$  given by $\langle S|T \rangle=S^*T$.

 The first question that requires settling is to answer  whether every product system arises this way which will then reduce the study of $E_0$-semigroups to that of product systems.  We have the following to say concerning the current status of the problem.
\begin{theorem}[\cite{Sundar_Existence}]
\label{existence}
Let $G$ be a locally compact group, and let $P \subset G$ be a subsemigroup with non-empty interior. Suppose that one of the following conditions hold:
\begin{enumerate}
\item $G$ is discrete, and $P$ is normal in $G$, i.e. $gPg^{-1}=P$ for all $g \in G$.
\item The group $G$ is abelian.
\item There exists $a \in P$ such that $\bigcup_{n=1}^{\infty}Pa^{-n}=G$.
\end{enumerate}
Then, the association $\alpha \to E_\alpha$ is a bijection between the class of $E_0$-semigroups (up to cocycle conjugacy) and the class of product systems over $P$ (up to isomorphism).
\end{theorem}
The above result was first proved by Arveson for $P=(0,\infty)$, and by the author and Murugan (\cite{Murugan_Sundar_continuous}) for the case of higher dimensional cones in $\mathbb{R}^d$.   We will explain the remaining results in the language of product systems. For, we can then translate them to $E_0$-semigroups in light of Thm. \ref{existence}.

Let $P \subset \mathbb{R}^d$ be a closed convex cone in $\mathbb{R}^d$.  For an isometric representation $V$ of $P$, let $E^V(x)$ and $F^V(x)$ be the symmetric and the antisymmetric Fock space of $Ker(V_x^*)$. Then, $E^V:=\{E^V(x)\}_{x \in P}$ is a product system, where the product rule on the exponential vectors is given by $e(\xi) \cdot e(\eta):=e(\xi+V_x\eta)$ for $\xi \in Ker(V_x^*), \eta \in Ker(V_y^*)$. We can make $F^V$ a product system  by defining a product as follows:
\[
(\xi_1 \wedge \xi_2 \wedge \cdots \wedge \xi_m) \odot (\eta_1 \wedge \eta_2 \wedge \cdots \wedge \eta_n)=V_{x}\eta_1 \wedge V_x\eta_2 \wedge \cdots \wedge V_x\eta_n \wedge \xi_1 \wedge \xi_2 \wedge \cdots \wedge \xi_m.
\]
The product systems $E^V$ and $F^V$ are the simplest examples of product systems, and they are  called CCR (canonical commutation relation) and CAR (canonical anticommutation relation) flows respectively.

A deep classification result of Arveson (\cite{Arv_path}), in the $1$-parameter case, states that decomposable product systems\footnote{ This means that  each fibre $E(x)$ has enough decomposable vectors. A vector $u \in E(x)$ is decomposable if $x=y+z$, then there exists $v \in E(y)$ and $w \in E(z)$ such that $u=vw$.} are precisely those of CCR flows. This is not true in the higher dimensional case, and we have the following structure theorem.

\begin{theorem}[\cite{Namitha_Sundar}, \cite{Injectivity}]
\label{Namitha}
Let $\Gamma:P \times P \to H$ be a weakly measurable map such that $\Gamma(x,y) \in Ker(V_y^*)$ for every $x,y \in P$. Suppose that there exists a measurable map $\alpha:P \times P \to \mathbb{T}$ such that for $x,y,z \in P$,
\begin{align*}
\Gamma(x,y+z)+V_y\Gamma(y,z)&=\Gamma(x,y)+V_y\Gamma(x+y,z), \textrm{~and~}\\
\frac{\alpha(x,y)\alpha(x+y,z)}{\alpha(x,y+z)\alpha(y,z)}&=e^{i Im \langle \Gamma(x,y+z)|V_y\Gamma(y,z)\rangle}.
\end{align*}
Define a new product $\boxdot $ on $E^V$ by $u \boxdot v=\alpha(x,y)W(V_x\Gamma(x,y))(u \cdot v)$, where $W(\cdot)$ are the  Weyl operators. Then, $(E^V,\boxdot)$ is a decomposable product system, and every decomposable product system (up to isomorphism) is of this form. If a decomposable product system has a unit\footnote{A unit is a non-zero cross section that is multiplicative.}, then it is isomorphic to the product system of a  CCR flow.
\end{theorem}
In the $1$-parameter case, up to isomorphism, we can take $\Gamma=0$.  This is not true  in the higher dimensional case as was shown  in \cite{Namitha_Sundar} by an  explicit computation of  the space of $2$-cocycles $\Gamma$ (modulo coboundaries) for a class of examples.

 Other  results, which are of equal importance as Thm. \ref{Namitha}, in the multiparameter theory are listed below; some are in total contrast to the $1$-parameter case.
\begin{enumerate}
\item The `map' $V \to E^V$ is injective (\cite{Injectivity}, \cite{Piyasa_Sundar}), and so is the `map' $V \to F^V$ (\cite{Namitha_CJM}).
\item There are uncountably many mutually inequivalent isometric representations $V$ of $P$ such that
\begin{itemize}
\item $E^V$ is not type one\footnote{This means that $E^V$ does not have enough units. }(\cite{Anbu}, \cite{Anbu_Sundar}),
\item $E^V$ is of type one, prime and has index $k$ for any given $k \in \mathbb{N}_0 \cup \{\infty\}$ (\cite{Piyasa_Sundar}, \cite{Piyasa_NYJM}),
\item  $F^V$ is not isomorphic to $E^V$ (\cite{Anbu_Decomposable1}, \cite{Namitha_CJM}, \cite{Vasanth}), and
\item the gauge group of $F^V$ does not act transitively on the set of units.
\end{itemize}
\item The isometric representations of $P$ with commuting range projections  such that $F^V$ is type I are precisely those that are direct sums of pullbacks of the standard shift semigroup $\{S_t\}_{t \geq 0}$ on $L^2([0,\infty))$. In particular,  in contrast to CCR flows, only pullbacks of the shift semigroup can give rise to type I, index one CAR flows if we restrict attention to semigroups of isometries with commuting range projections (\cite{Namitha_CJM}).
\end{enumerate}

\begin{center}
\textsf{Future directions}
\end{center}

The above results indicate that the multiparameter theory is vastly different from that of the $1$-parameter case. This makes it a fertile ground for future research. The author wishes to  pose a couple of concrete questions in this regard.
\begin{enumerate}
\item Is it true that  the conclusion of Thm. \ref{existence} holds for every  subsemigroup $P$  of a locally compact group?
\item Is there an example of an isometric representation $V$ of a cone $P$ (say $\mathbb{R}_{+}^{2}$) such that $V$ is not a pullback of the shift semigroup but  $F^V$ is type one and has index one?
\end{enumerate}

 Power's  approach (\cite{Powers_TypeIII}, \cite{Powers_Index}, \cite{Powers}, \cite{Powers_CPflow}) to $E_0$-semigroups, and the probabilistic approach of Tsireslon  ( \cite{Tsi}, \cite{Tsirelson}, \cite{Tsirelson2003}) via random sets, stochastic processes, point processes etc. to product systems remain unexplored in higher dimensions. The author is fairly confident that the analysis of the resulting models of $E_0$-semigroups and product systems will  result in deep and interesting mathematics.

\section{$C^*$-algebras, $K$-theory and Dynamics}

\noindent \textit{By Prahlad Vaidyanathan [prahlad@iiserb.ac.in] from IISER Bhopal.}

\bigskip

A $C^*$-algebra is a norm-closed self-adjoint subalgebra of the space of bounded linear operators on a complex Hilbert space. These algebras were introduced in the 1930s by von Neumann, Gelfand, and many others. In the 1970s, the subject received a dramatic revitalization due to the work of Brown-Douglas-Fillmore, Elliott, Kasparov and others. Specifically, the introduction of $K$-theory and related tools from algebraic topology have brought two disparate areas closer together, and deepened our understanding of both. In what follows, we will describe the work that we, along with our collaborators, have done in this field in the recent past.

\textbf{Nonstable $K$-theory:} $K$-theory is a pair of abelian groups associated to a $C^*$-algebra. These functors are powerful invariants and are central to the classification programme. However, in constructing these groups, one necessarily loses a great deal of information. Specifically, the $K$-groups are inductive limits of certain homotopy groups associated to the $C^*$-algebra \cite{thomsen1991}. Nonstable $K$-theory is the branch of study that attempts to understand these homotopy groups by themselves, and how they behave before one passes to the limit.

The subject began in the 1980s, when Rieffel \cite{rieffel1983} introduced the notions of connected and general stable rank, which control the behaviour of the maps in the inductive limit mentioned above. Nica \cite{nica2011} further developed the theory of these two stable ranks, and proved that they are homotopy invariant. In \cite{vaidyanathan2019a}, we determined how these ranks behave with respect to some natural constructions, most notably pullbacks. As a result, we were able to compute these ranks for various classes of $C^*$-algebras. In \cite{nirbhay2021a}, we estimated these ranks for $C(X)$-algebras, and for a variety of group $C^*$-algebras and crossed products.

While computing these ranks helps us understand the behaviour of nonstable $K$-groups, we do not have many tools to compute these groups explicitly. Therefore, we began studying $K$-stability; a phenomenon exhibited by a number of infinite dimensional $C^*$-algebras, where the nonstable $K$-groups naturally coincide with the stable $K$-groups.

In \cite{seth2020}, we showed that a continuous $C(X)$-algebra with K-stable fibers must itself be K-stable. In \cite{seth}, we showed that AF-algebras are K-stable provided they have slow dimension growth. We were also able to compute the rational homotopy groups of the unitary groups of finite dimensional algebras. Using that, we showed that an AF-algebra is $K$-stable if and only if it is rationally K-stable; a condition that is much easier to verify. In \cite{seth2023a}, we showed that A$\mathbb{T}$-algebras exhibit the same phenomenon, raising an interesting question about AH-algebras in general. In an attempt to understand this further, we showed in \cite{seth2023} that the notions of rational K-stability and $K$-stability are indeed distinct in general, and that the results of \cite{seth2020} were also true for rational $K$-stability.

At the moment, the subject has a number of open questions. We do not, as yet, have explicit calculations of nonstable $K$-groups, nor do we know how these groups behave under natural constructions. Moreover, we do not know if all simple, real rank zero $C^*$-algebras are $K$-stable. This is a long-standing conjecture of Zhang \cite{zhang1991}.

\textbf{Rokhlin Dimension:} The study of group actions on $C^*$-algebras is an integral part of operator algebra theory. Given a group action on a $C^*$-algebra, one is often interested in the structure of the associated crossed product $C^*$-algebra, and an important question in this context is to determine when certain ‘regularity’ properties pass from the underlying algebra to the crossed product. These include many properties that are useful from the point of view of the classification programme such as finiteness of nuclear dimension, simplicity, $\mathcal{Z}$-stability, etc. Rokhlin dimension, introduced in \cite{hirshberg2015}, is a notion analogous to covering dimension for group actions which has proved useful in this context.

This has been studied for actions of finite groups and the integers \cite{hirshberg2015}, compact groups \cite{gardella2017} and for residually finite groups \cite{szabo2019}. In each case, it was proved that actions with finite Rokhlin dimension allow us to deduce a number of regularity properties of the crossed product $C^*$-algebra from those of the underlying algebra (particularly in the compact case, and under the ‘commuting towers’ hypothesis).

In \cite{vaidyanathan2022}, we showed that $C(X)$-linear actions of compact groups on $C(X)$-algebras have finite Rokhlin dimension provided the induced action on their fibers have finite Rokhlin dimension. This allowed us to produce new classes of classifiable $C^*$-algebras associated to equivariant vector bundles. In \cite{m2024}, we studied Rokhlin dimension for actions of residually finite groups, building on the work of \cite{szabo2019}. We showed a number of useful permanence properties, and also described the ideal structure of the associated crossed product. In \cite{m2024a}, we studied automorphisms of $AF$-algebras, and showed that under certain conditions, the corresponding crossed product is an A$\mathbb{T}$-algebra.

Once again, there are a number of interesting avenues of research in this field. Even in the well-studied case of integer actions, it is unclear what values of Rokhlin dimension are admissible. This is primarily because we do not have a computable obstruction which would provide us with lower bounds for the Rokhlin dimension. It would also be interesting to look at certain natural actions on specific classes of $C^*$-algebras, and compute their Rokhlin dimension. Finally, we have yet to give an appropriate definition of Rokhlin dimension for actions of locally compact groups. These are amongst the questions that we plan to pursue in the future.

In conclusion, the projects described above involve a synthesis of ideas from topology, algebra and analysis, and therefore have a wide appeal. A number of new tools have been discovered in these fields in the past decade, and this has led to renewed interest in these problems. Nonstable $K$-theory is a mature field, and the open problems require genuinely novel ideas. Rokhlin dimension, on the other hand, is relatively young and exciting. It is hoped that many of those outstanding problems will be resolved in the coming years.

\end{document}